\documentclass[a4paper,11pt]{article}
\usepackage{amsmath,amsthm,amssymb}
\usepackage[mathscr]{eucal}
\theoremstyle{plain}
\newtheorem*{thm*}{Theorem}
\newtheorem{lem}{Lemma}
\newtheorem*{stats*}{Statements}
\theoremstyle{definition}
\newtheorem*{Def*}{Definition}
\newtheorem*{defA*}{Condition A}
\newtheorem*{defB*}{Condition B}
\newtheorem*{defC*}{Condition C}
\newtheorem*{defS*}{Scattering data}
\newtheorem*{nots*}{Notations}
\newtheorem*{rem*}{Remark}
\newtheorem*{rems*}{Remarks}

\providecommand{\BS}[1]{\boldsymbol{#1}}
\providecommand{\C}[1]{\mathcal{#1}}
\providecommand{\D}[1]{\mathbb{#1}}
\providecommand{\E}[1]{\mathscr{#1}}
\providecommand{\R}[1]{\mathrm{#1}}
\newcommand{\DS}{\displaystyle}
\newcommand{\dd}{\mathrm{d}}
\newcommand{\ee}{\mathrm{e}}
\newcommand{\ii}{\mathrm{i}}
\renewcommand{\Im}{\operatorname{Im}}
\newcommand{\ord}{\mathrm{O}}
\renewcommand{\Re}{\operatorname{Re}}
\DeclareMathOperator{\res}{Res}
\begin{document}
\title{{\Large\textbf{Characteristic properties of the scattering data\\
                       for the mKdV equation on the half-line}}}
\author{{\normalsize Anne \textsc{Boutet} de \textsc{Monvel}$^\ast$
and Vladimir
\textsc{Kotlyarov}$^\dagger$}\\[1mm]
{\scriptsize $^{\ast}$ Institut de Math\'ematiques de Jussieu, case 7012, 
     Universit\'e Paris 7,}\\[-2mm] 
{\scriptsize 2 place Jussieu, 75251 Paris, France}\\
{\scriptsize 
$^{\dagger}$ Mathematical Division, Institute for Low Temperature Physics,}\\[-2mm]
{\scriptsize 47 Lenin Avenue, 61103 Kharkiv, Ukraine}}
\date{}
\maketitle
\begin{abstract}
In this paper we describe characteristic properties of the
scattering data of the compatible eigenvalue problem for the pair of differential
equations related to the modified Korteweg-de Vries (mKdV) equation
whose solution is defined in some
half-strip $(0<x<\infty)\times[0,T]$, or in the quarter plane
$(0<x<\infty)\times(0<t<\infty)$.
We suppose that this solution has a $C^{\infty}$ initial
function vanishing as $x\to\infty$, and $C^{\infty}$ boundary
values, vanishing as $t\to\infty$ when $T=\infty$.
We study the corresponding scattering problem for the compatible
Zakharov-Shabat system of differential equations associated with the mKdV
equation and obtain a representation of the solution of the mKdV equation
through Marchenko integral equations of the inverse scattering method. 
The kernel of these
equations is valid only for $x\geq 0$ and it takes into account all specific properties
of the pair of compatible differential equations in the chosen half-strip or
in the quarter plane. The main result of the paper is the collection
A\textendash B\textendash C of characteristic properties of the scattering
functions given below.
\end{abstract}
\section{Introduction}
\subsection{}
Initial value problems on the whole line for nonlinear integrable equations
such as the nonlinear Schr\"odinger equation,
the Korteweg-de Vries equation, the sine-Gordon equation, etc.\ are
well studied. The solvability of the Cauchy problem, multi-soliton solutions,
the proof that these nonlinear equations are completely integrable
infinite-dimensional Hamiltonian systems are the most significant results in
the soliton theory on the whole line.

At the same time the initial-boundary value problem on the half-line for
nonlinear integrable equations has not been studied so far. In
the last decade attention to those problem has strongly increased.
Among papers \cite{AS1}-\cite{DMS1}, \cite{F}-\cite{H2},
\cite{Sab99}-\cite{T} devoted to
this problem, the most interesting  results were obtained by
A.S.~Fokas \cite{F},  A.S.~Fokas and A.R.~Its \cite{FI2}-\cite{FI}.
Later, in \cite{F97} A.S.~Fokas has proposed a general method for solving
boundary value problems for two-dimensional linear and integrable nonlinear
partial differential equations.
This method, which was  further developed in \cite{Fo}, 
\cite{F011}-\cite{F01}, is based on the simultaneous spectral
analysis of the two eigenvalue equations of the associated Lax pair.
It expresses  the solution in terms of the solution of a matrix
Riemann-Hilbert  problem in the complex plane of the
spectral parameter. The spectral functions determining the Riemann-Hilbert
problem are expressed in terms of the initial and  boundary values of the
solution. The fact that these initial and boundary
values are in general related can be expressed in a simple way in terms
of a global relation satisfied by the corresponding spectral functions.

In the framework of this approach we recently
found characteristic properties of the scattering data for the compatible
Zakharov-Shabat eigenvalue problem associated with focusing and defocusing
nonlinear Schr\"{o}dinger equations on the half-line with initial and
boundary functions of Schwartz type \cite{BMK}.

Recently in \cite{BFS} (see also \cite{H2}) an initial-boundary value
problem for the mKdV equation on the half-line was analyzed by expressing
the solution in terms of the solution of a matrix Riemann-Hilbert problem
in the complex $k$-plane. In particular,
it is shown that for a subclass of boundary conditions,
the ``linearizable boundary conditions'', all spectral functions can be
computed from the given initial data by using  algebraic manipulations of some
``global relation''. Thus in this case, the problem on the half-line can be solved
as efficiently as the problem on the whole line. 

But the general
initial-boundary value problem on the half-line remains non-linea\-rizable.
Characteristic properties of the spectral functions were not considered in
\cite{BFS}. In this connection the characterization of the spectral
functions becomes important. Besides, a description of the characteristic
properties of the scattering or spectral data is a very important problem
in itself \cite{Mar}.

Most of papers initiated of problems on the half-line deal with nonlinear
dynamics of the spectral or scattering data. We prefer the approach of
A.S.~Fokas and A.R.~Its where scattering (spectral) data have  trivial dynamics. 
But then analytic properties of the
scattering data are more complicated. Therefore it is necessary
to give their complete description when, of course, initial and boundary
functions belong to suitable classes of functions. 
We introduce spectral
data in a natural way as in ordinary scattering problems:
\begin{enumerate}
\item
First, a ``scattering
matrix" for the $x$-equation (by initial function).
\item
Then, a  ``scattering
matrix" for the $t$-equation (by boundary functions),
\item
Finally, a ``scattering matrix" for the compatible $x$ and $t$-equations 
as by-product.
\end{enumerate}
In this case  a kernel
of the Marchenko integral equations or a jump matrix of the corresponding
Riemann-Hilbert problem has an explicit $x$,  $t$ dependence.
That makes possible to study the asymptotic behavior of the solution of the
non-linear problem by using, for example, the powerful steepest descent method
of  P.~Deift and X.~Zhou \cite{DZ1,DZ2} while the non-linear dynamics of the
spectral or scattering data makes almost impossible to obtain an effective
asymptotics of the solution.

\subsection{}

In this paper we consider the problem to characterize ``scattering
data'' for a compatible pair of 
differential equations attached to the modified Korteweg-de Vries (mKdV)
equation.
Let $q(x,t)$ be a real-valued solution 
of the mKdV equation
\begin{align}                                  \label{mkd}
&q_t+q_{xxx}-6\lambda q^2q_x=0,\\
&x\in\D{R}_+, \quad
t\in [0,T], \  T\leq\infty, \quad \lambda=\pm 1\notag
\end{align}
in the half-strip or quarter $xt$-plane and suppose the initial function
\[
q(x,0)=u(x)\text{ with } x\in\D{R}_+;
\]
and the boundary values
\[
q(0,t)=v(t)\quad
q_x(0,t)=v_1(t) \quad
q_{xx}(0,t)=v_2(t)
\text{ with }t\in[0,T], \quad T\leq\infty
\]
are $C^{\infty}$, and $u(x)\in\C{S}(\D{R}_+)$, 
where $\C{S}(\D{R}_+)$ is 
the Schwartz space of rapidly decreasing functions on
$\D{R}_+$, i.e.\ $C^{\infty}$ functions whose derivatives of any order $n\geq 0$
vanish at infinity faster than any negative power of $x$.
For $T=\infty$ the boundary values are also supposed to be rapidly
decreasing: $v(t), v_1(t), v_2(t)\in\C{S}(\D{R}_+)$.

\begin{rem*}
We consider here the IBV problem for the mKdV
equation (\ref{mkd}) in the first quarter $(x\geq 0,\, t\geq 0)$ of the
$xt$-plane. This problem differs from that studied in \cite{BFS}
which is also on the first quarter but for the mKdV equation of the form:
\[
q_t-q_{xxx}+6\lambda q^2q_x=0.
\]
This form can be easily reduced to (\ref{mkd}), but then the IBV problem
is on the second quarter  $(x\leq 0,\, t\geq 0)$ of the $xt$-plane. 
Scattering (spectral) data for that
problem have different analytic properties. 
It is well-known that for KdV
and mKdV equations there are
differences between the IBV problems for $x>0$ and for $x<0$.
\end{rem*}

To study the solution $q(x,t)$ we shall use  spectral
analysis of a compatible eigenvalue problem for the linear $x$-equation
\begin{align}  \label{xeq}
& w_x+\ii k\sigma_3w=Q(x,t)w, \\
&\sigma_3 =
\begin{pmatrix}
1 & 0 \\
0 & -1
\end{pmatrix}
,  \notag \\
&Q(x,t) =
\begin{pmatrix}
0 & q(x,t) \\
\lambda q(x,t) & 0
\end{pmatrix}
\notag
\end{align}
and for the linear $t$-equation
\begin{align}  \label{teq}
&w_t+4\ii k^3\sigma_3w=\hat Q(x,t,k)w, \\
&\hat Q(x,t,k)=2Q^3(x,t)-Q_{xx}-2\ii k(Q^2(x,t)+Q_x(x,t))\sigma_3+4k^2Q(x,t). 
\notag
\end{align}
This is the well-known Ablowitz-Kaup-Newel-Segur \cite{AS} or
Zakharov-Shabat \cite{ZS}  system of linear equations, which are
compatible if and only if $q(x,t)$ satisfies the mKdV
equation.

The main goal of the present paper is to study the scattering problem
for compatible differential equations (\ref{xeq}) and (\ref{teq})
on the half-strip or on the quarter of the $xt$-plane. We will combine
the Marchenko integral equation and the corresponding Riemann\textendash Hilbert problem
in our approach to obtain characteristic properties of scattering data.
We get a description of these characteristic properties and obtain a
representation of the solution of the mKdV equation through
Marchenko integral equations of the inverse scattering method.  The kernel
of these equations is valid for $x\geq 0$ only and  it
takes into account all specific properties occurred for compatible differential
equations.  In particular,
the solution
of the mKdV equation given by our method does not have continuation for $x<0$.
It is well defined on the half-strip or on the quarter of the $xt$-plane.
Then, if one is using our representation of
the solution of the mKdV equation and the explicit
$(x,t)$-dependence of the kernel of the Marchenko integral equations or the
jump matrix in the Riemann-Hilbert problem one can easily obtain the asymptotic
behavior of the solution in the same way as in \cite{DZ1,DZ2} or
\cite{K1}.

\subsection{Scattering data: definition
            and properties}

Let $q(x,t)$ be a real-valued solution of Equation (\ref{mkd})
with initial and boundary functions satisfying smoothness and decreasing
assumptions described above.
Let us define $\Sigma=\{k\in\D{C}\mid\Im k^3=0\}$ and domains $\Omega_1,\dots,
\Omega_6$ as depicted on Figure 1.
\begin{figure}[ht]
\setlength{\unitlength}{1mm}
\begin{picture}(120,47)
\thicklines
\put(70,22){\line(-2,-3){16}}
\put(70,22){\line(2,3){16}}
\put(70,22){\line(-2,3){16}}
\put(70,22){\line(2,-3){16}}
\put(38,22){\line(1,0){64}}
\put(70,22){\line(-1,0){10}}
\put(85,30){$\Omega_1$}
\put(85,12){$\Omega_6$}
\put(50,30){$\Omega_3$}
\put(68.5,39){$\Omega_2$}
\put(50,12){$\Omega_4$}
\put(68.5,5){$\Omega_5$}
\end{picture}
\caption{$\Sigma=\{k\in\D{C}\mid\Im k^3=0\}$ and $\Omega_1,\dots,\Omega_6$.}
\label{fig.mkd-1}
\end{figure}
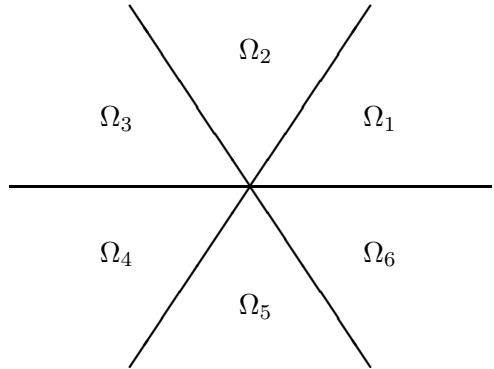

\noindent
Scattering data are introduced as follows.

1.
The initial function $u(x)=q(x,0)$ and the
$x$-equation (\ref{xeq}) with
$t=0$ define the Jost solution $\Psi(x,0,k)=\exp(-\ii kx\sigma_3)
+\R{o}(1)$, $x\to\infty$  and  a ``scattering matrix"
\[
S(k):=\Psi^{-1}(0,0,k)=\begin{pmatrix} s_2^+(k)&-s_1^+(k)\\-s_2^-(k)&s_1^-(k)
\end{pmatrix}, \quad s_1^-(k)=\bar s_2^+(\bar k),
\quad s_2^-(k)=\lambda\bar s_1^+(\bar k).
\]
In particular, they define
\begin{itemize}
\item
the spectral function 
$r(k)=-s_2^-(k)/s_2^+(k)$, called the ``reflection coefficient'',
\item
eigenvalues $k_j\in\D{C}_+$, $j=1,\dots,n$, which are zeros of $s_2^+(k)$,
\item
numbers $m_j=[\ii s_1^+(k_j)\dot s_2^+(k_j)]^{-1}$,
$k_j\in\D{C}_+$.
\end{itemize}

2.
Boundary data $v(t)= q(0,t)$, $v_1(t)=q_x(0,t)$, $v_2(t)=q_{xx}(0,t)$,
with $t\in[0,T]$,  $T\leq\infty$ and the $t$-equation with $x=0$
define a solution $Y(0,t,k)=\exp(-4\ii k^3t\sigma_3)$, $t\geq T$, then
a ``scattering matrix"
\[
P(k):=Y(0,0,k)=
\begin{pmatrix}   
p_1^-(k)& p_1^+(k)\\
p_2^-(k)&p_2^+(k),
\end{pmatrix}
\quad p_1^-(k)=\bar p_2^+(\bar k),
\quad p_2^-(k)=\lambda\bar p_1^+(\bar k),
\]
and, together with $S(k)$, another
``scattering matrix"
\[
R(k)=S(k)P(k) = 
\begin{pmatrix}   
r_1^-(k)& r_1^+(k)\\r_2^-(k)&r_2^+(k)
\end{pmatrix} .
\]
Now we introduce:
\begin{itemize}
\item
one more spectral function
$c(k)=\dfrac{p_2^-(k)}{s_2^+(k)r_1^-(k)}$, $k\in\Omega_2$,
\item
eigenvalues $z_j\in\Omega_2$, $j=1,\dots,m$, which are zeros of $r_1^-(k)$,
\item
numbers
$m_j^2=-\ii\res_{k=z_j}c(k)$ ($z_j\in\Omega_2$),
which depend on the initial and boundary functions.
\end{itemize}

\begin{defS*}
We define the set
\[
\C{R}=\{k_1,\dots,k_n\in\D{C}_+;\ z_1,\dots, z_m\in\Omega_2;\
r(k),\;  k\in\D{R};\ c(k),\;  k\in\Omega_2\}
\]
as ``scattering data'' of the compatible eigenvalue problem for the system of
differential equations (\ref{xeq})-(\ref{teq}) with $q(x,t)$ satisfying the
mKdV equation
(\ref{mkd}).
\end{defS*}

\subsubsection*{Recovering of 
    $\BS{q(x,t)}$ from scattering data}

Then we prove that the solution $q(x,t)$ of the non-linear
problem (\ref{mkd}) can be written
\begin{equation}
q(x,t)=-2\lambda K_2(x,x,t)  \label{q=K2}
\end{equation}
where $K_2(x,y,t)$, together with $K_1(x,y,t)$, satisfies the
Marchenko integral equations:
\begin{align}
& K_1(x,y,t)+\lambda\int_{x}^{\infty}{K}_2(x,z,t)H(z+y,t)\dd z=0
\text{ for }0\leq x<y<\infty ,  \label{ME2} \\
& K_2(x,y,t)+H(x+y,t)+\int_{x}^{\infty}K_1(x,z,t)H(z+y,t)\dd z=0
\end{align}
with kernel
\begin{align}	                                          \label{H0}
H(x,t)&=\frac{1-\lambda}{2}\left(
\sum_{\substack{k_j\in\Omega_1\cup\Omega_3}}
m_j\ee^{\ii k_jx+8\ii k_j^3t}
+\sum_{\substack{z_j\in\Omega_2}}
m_j^2\ee^{\ii z_jx+8\ii z_j^3t}\right)
\notag\\
&\quad
+\frac{1}{2\pi}\int_{\partial\Omega_2}
c(k)\ee^{\ii kx+8\ii k^3t}\dd k
\quad
+\frac{1}{2\pi}\int_{-\infty}^{\infty}r(k)\ee^{\ii kx+8\ii k^3t}\dd k.
\end{align}

\subsubsection*{Properties of 
                scattering data}
Now we introduce three sets of conditions on a set 
$\C{R}$
of numbers 
$k_1,\dots,k_n\in\D{C}_+$, $z_1,\dots,z_m\in\Omega_2$ and functions
$r(k)$, $k\in\D{R}$, $c(k)$,  $k\in\Omega_2$.

\begin{defA*}
\textit{Conditions on} $r(k)$, $k\in\D{R}$:
\begin{itemize}
\item  
$r(k)\in C^{\infty}(\D{R})$, $r(-k)=\bar r(k)$;
\ $r(k)=\ord(k^{-1})$,\ as
$k\to\infty$;\  $|r(k)|<1$, if $\lambda=1$.
\item  
$r(k)=\dfrac{-s_2^-(k)}{s_2^+(k)}$, where $s_2^-(k)$ is
analytic in $k\in\D{C}_-$, $s_2^+(k)$ is analytic for 
$k\in\D{C}_+$ and has the form:
\[
s_2^+(k)=
\biggl(\prod_{j=1}^{n}
\frac{k-k_{j}}{k-\bar{k}_{j}}\biggr)^{\frac{1-\lambda}{2}}
\exp {\biggl[ \frac{\ii}{2\pi}
\int_{-\infty}^{\infty}\frac{\log (1-\lambda|r(\mu)|^{2})
\dd\mu}{\mu-k}\biggr]},\quad k\in\D{C}_+.
\]
\item 
The function given by
\[
\int_{-\infty}^{\infty}r(k)\;\ee^{\ii k(x+8k^2t)}\dd k
\]
is $C^{\infty}$ in $x,t$ for $x>0$
and $t\geq 0$.
\end{itemize}
\end{defA*}

\begin{defB*}
\textit{Conditions on} 
$\E{K}=\{k_1,\dots,k_n\in\D{C}_+;\ z_1,\dots,z_m\in\Omega_2\}$:
\begin{itemize}
\item 
If $\lambda=1$, then $\E{K}=\varnothing$.
\item 
If $\lambda=-1$, then $\E{K}$ satisfies symmetry conditions:
\begin{alignat*}{4}
&k_j=\ii\kappa_j,&\quad&1\leq j\leq n_1\leq n=n_1+2n_2
&\qquad\;&k_{n_1+l}=-\bar k_{n_1+n_2+l},&\quad&1\leq l\leq n_2\\
&z_j=\ii\mu_j,&&1\leq j\leq m_1\leq m=m+2m_2,
&&z_{m_1+l}=-\bar z_{m_1+m_2+l},&&1\leq l\leq m_2.
\end{alignat*}
\end{itemize}
\end{defB*}

\begin{defC*}
\textit{Conditions on} $c(k)$, $k\in\Omega_2$:
\begin{itemize}
\item 
If $\lambda=1$, then $c(k)$ is analytic in $k\in\D{C}_+$ ($T<\infty$) or in
$\Omega_2$ ($T=\infty$) and it is bounded on $\overline\Omega_2$.
\item 
If $\lambda=-1$, then $c(k)$ is meromorphic
in the half-plane $\D{C}_+$ ($T<\infty$)
or in $\Omega_2$ ($T=\infty$), where it
has poles at $z_1,z_2,\dots,z_m$.
\item  
$c(k)=-\bar c(-\bar k)$, and
$c(k)\to 0$, $k\to\infty$, $k\in\overline{\D{C}}_+$ or
$k\in\overline\Omega_2$.
\item  
$c(k)$ has $C^{\infty}$ boundary values on $\D{R}$ or $\partial\Omega_2$.
\item 
If $T=\infty$, then $\dfrac{\dd^n c(k)}{\dd k^n}\bigg\vert_{k=0}=
-\dfrac{\dd^n r(k)}{\dd k^n}\bigg\vert_{k=0}$, $n=0,1,2,\ldots$.
\item 
The function given by
\[
\int_{\partial\Omega_2}
c(k)\ee^{8\ii k^3t}\dd k +\int_{-\infty}^{\infty}r(k)\ee^{8\ii k^3t}\dd k
\]
is $C^{\infty}$ in $t$.
\end{itemize}
\end{defC*}

\subsection{Main theorem}

The main result is that properties A\textendash B\textendash C are characteristic:

\begin{thm*} 
Conditions \emph{A\textendash B\textendash C} on 
\[
\C{R}=\{k_1,\dots,k_n\in\D{C}_+;\ z_1,\dots, z_m\in\Omega_2;\
r(k),\,  k\in\D{R};\ c(k),\,  k\in\Omega_2\}
\]
characterize the scattering data of the
compatible eigenvalue problem \emph{(\ref{xeq})-(\ref{teq})} for $x$- and 
$t$-equations defined by
a solution $q(x,t)$ of the mKdV equation $(\ref{mkd})$ with initial
function $u(x)\in\C{S}(\D{R}_+)$ and boundary values
$v(t), v_1(t), v_2(t)\in C^{\infty}[0,T]$ if $T<\infty$, or
$v(t), v_1(t), v_2(t)\in\C{S}(\D{R}_+)$ if $T=\infty$.
\end{thm*}

\begin{rem*}
The third item in condition A (about smoothness of the
integral with respect to $x$ and $t$) is fulfilled if the initial and
boundary functions obey the following relations~:
\begin{equation}  \label{dn0}
\frac{\dd^n}{\dd x^n}u(x)\!\Bigm\vert_{x=0}= 
\frac{\dd^n}{\dd t^n}v(t)\!\Bigm\vert_{t=0}= \frac{\dd^n}{\dd t^n}
v_1(t)\!\Bigm\vert_{t=0}= \frac{\dd^n}{\dd t^n}
v_2(t)\!\Bigm\vert_{t=0}= 0
\end{equation}
for any $n\geq 0$. In this case the reflection coefficient $r(k)$ is in
the Schwartz space $\C{S}(\D{R})$. Under these assumptions
the last item of condition C is also fulfilled.
\end{rem*}

\begin{rem*}
If $u(x)\equiv 0$, conditions A  are
trivial: $r(k)\equiv 0$, $a(k)\equiv 1$,
$\{k_1,\dots,k_n\}=\varnothing$. There
remain only conditions B and C. They mean that any solution
of the mKdV equation in the half-strip or quarter plane with zero initial
function is parametrized by only one function $c(k)=-\bar c(-\bar k)$,
analytic
($\lambda=1$) or meromorphic ($\lambda=-1$) in $\Omega_2$ 
with poles at $z_1,\dots,z_m$.
\end{rem*}

\section{Basic solutions of 
     compatible $x$- and $t$-equations}
\setcounter{equation}{0}

Let us write the $x$- and $t$-equations in the form
\begin{align}  \label{U}
W_x&=U(x,t,k)W, \\ \label{V}
W_t&=V(x,t,k)W,
\end{align}
where $U(x,t,k)$ and $V(x,t,k)$ are matrices given by
\begin{align*}
U(x,t,k)&=Q(x,t)-\ii k\sigma_3, \\
V(x,t,k)&=2Q^3(x,t)-Q_{xx}-2\ii k(Q^2(x,t)+Q_x(x,t))\sigma_3+4k^2Q(x,t)-
4\ii k^3\sigma_3.
\end{align*}

\begin{lem}                                         \label{lem.2.1} 
Let the system $(\ref{U})$, $(\ref{V})$ be compatible for all
$k$.
Let $W(x,t,k)$ satisfy the $x$-equation $(\ref{U})$ for all $t$, and let 
$W(x_0,t,k)$ satisfy the $t$-equation $(\ref{V})$ for some $x=x_0$
(including the case $x_0=\infty $). Then $W(x,t,k)$ satisfies the 
$t$-equation for all $x$.
\end{lem}

\begin{proof}
See e.g.~\cite{BMK}.
\end{proof}

\begin{nots*}
The over-bar denotes the
complex conjugation. $\D{C}_{\pm}$ denotes the upper (lower) complex half
plane. If $A=\begin{pmatrix}A^-&A^+\end{pmatrix}$ 
denotes a $2\times 2$
matrix, the vectors $A^{\mp}$ denote the first and second columns of $A$. 
We also denote
$[A,B]=AB-BA$.
\end{nots*}

In this section we shall introduce basic solutions of compatible
$x$- and $t$-equations. 

\subsection{First basic solution}

The first basic solution is a matrix-valued Jost
solution of the $x$-equation (\ref{xeq}).
It has the triangular integral representation
(see e.g.~\cite{FT})
\begin{equation}  \label{PSI}
\Psi(x,t,k)=
\left(\ee^{-\ii kx\sigma_3}+\int_x^{\infty}K(x,y,t) 
\ee^{-\ii ky\sigma_3}\dd y\right)\ee^{-4\ii k^3t\sigma_3},
\end{equation}
where real-valued matrix $K(x,y,t)$ has the form
\[
K(x,y,t)=
\begin{pmatrix}
K_1(x,y,t) & \lambda K_2(x,y,t) \\
K_2(x,y,t) & K_1(x,y,t)
\end{pmatrix}
\]
with entries in 
$C^{\infty}(\D{R}_+\times\D{R}_+\times\D{R}_+)$
and rapidly decreasing as 
$x+y\to\infty$ for any $t\in\D{R}_+$. 
Matrices $K(x,x,t)$ and $Q(x,t)$ are connected by the relation:
\begin{equation}  \label{sigma}
[\sigma_3, K(x,x,t)]=Q(x,t)\sigma_3.
\end{equation}
This last equality yields important formula (\ref{q=K2}) for the solution 
$q(x,t)$ of the modified Korteweg-de Vries  equation. The matrix $\Psi(x,t,k)$
satisfies the $x$-equation (\ref{xeq}) and it satisfies the $t$-equation 
(\ref{teq}) with $x=\infty$, because the matrix $\ee^{-\ii k(x+4k^2t)\sigma_3}$ 
is a solution of both equations (\ref{xeq}) and 
(\ref{teq}) with $Q(x,t)\equiv 0$. Lemma \ref{lem.2.1} implies that 
$\Psi(x,t,k)$ satisfies the $t$-equation for any $x\in\D{R}_+$ due to the
compatibility of the $x$- and $t$-equations.

The triangular integral representation (\ref{PSI}) and Lemma
\ref{lem.2.1} imply the following properties of the matrix-valued
Jost solution $\Psi(x,t,k)$ (cf.\ \cite{FT}):

\subsubsection*{Properties of the 
                  first basic solution}
\begin{enumerate}
\item  
$\Psi (x,t,k)$ satisfies the $x$- and $t$-equations 
(\ref{xeq})-(\ref{teq}).
\item  
$\Psi(x,t,k)=\Lambda\bar{\Psi}(x,t,k)\Lambda^{-1}$ for $k\in\D{R}$, 
$\Lambda=
\begin{pmatrix}
0 & 1 \\
\lambda & 0
\end{pmatrix}
$.
\item  
$\det\Psi(x,t,k)\equiv 1\text{ for }k\in\D{R}$.
\item  
$(x,t,k)\mapsto 
\Psi(x,t,k)\in C^{\infty}(\D{R}_+\times\D{R}_+\times\D{R})$
\item  
$\Psi^+(x,t,k)$ is analytic in $k\in\D{C}_+$,
$\Psi^-(x,t,k)$ is analytic in $k\in\D{C}_-$.
\item 
$\Psi^\pm(x,t,k)=\overline{\Psi^\pm(x,t,-\bar k)}$
\item  
For $k\to\infty $,
\begin{align*}
& \ee^{\ii kx+4\ii k^{3}t}\Psi ^{-}(x,t,k)=
\begin{pmatrix}
1 \\
0
\end{pmatrix}
+\ord(k^{-1})\quad\text{if }\Im k\leq 0, \\
& \ee^{-\ii kx-4\ii k^{3}t}\Psi^+(x,t,k)=
\begin{pmatrix}
0 \\
1
\end{pmatrix}
+\ord(k^{-1})\quad \text{if }\Im k\geq 0.
\end{align*}
\end{enumerate}

\subsection{Second basic solution}

Now let us introduce the second basic solution $\Phi(x,t,k)$ of the $x$- and
$t$-equations which satisfies the initial condition
\begin{equation}  \label{PHI=E}
\Phi(0,0,k)=\sigma_0\equiv
\begin{pmatrix}
1 & 0 \\
0 & 1
\end{pmatrix}.
\end{equation}
It can be represented as a product of two matrices:
\begin{equation}  \label{PHI}
\Phi(x,t,k)=\varphi(x,t,k)\hat\varphi(t,k),
\end{equation}
where $\varphi(x,t,k)$ satisfies the $x$-equation under the condition 
$\varphi(0,t,k)=\sigma_0$, and $\hat\varphi(t,k)$ satisfies the $t$-equation
with $x=0$ under initial condition $\hat\varphi(0,k)=\sigma_0$. Lemma 
\ref{lem.2.1} implies that $\Phi(x,t,k)$ is a compatible solution of the $x$-
and $t$-equations. The existence of the solution $\varphi(x,t,k)$ and its
representation
\begin{equation}  \label{phi}
\varphi(x,t,k)=\ee^{-\ii kx\sigma_3}+\int_{-x}^x A(x,y,t)
\ee^{-\ii ky\sigma_3}\dd y
\end{equation}
by some integral kernel $A(x,y,t)$ are proved in \cite{BMK}. The matrix
$\hat\varphi(t,k)$ can be found as solution of the Volterra integral
equation:
\begin{equation}  \label{hatphi}
\hat\varphi(t,k)=\ee^{-4\ii k^3t\sigma_3}+\int_0^t
\ee^{4\ii k^3(\tau-t)}\hat Q(0,\tau,k)\hat\varphi(\tau,k)\dd\tau,
\end{equation}
where
\[
\hat Q(0,t,k)=
\left(\begin{smallmatrix}
-2\ii\lambda k v^2(t) &2\lambda v^3(t)+4k^2 v(t)+2\ii kv_1(t)-v_2(t) \\
2v^3(t)+4\lambda k^2 v(t)-2\ii\lambda k v_1(t)-\lambda v_2(t)&2\ii\lambda
k v^2(t).
\end{smallmatrix}\right)
\]
Besides $\hat\varphi(t,k)$ has the integral representation:
\begin{align}  \label{hphi}
\hat\varphi(t,k)
&=\ee^{-4\ii k^3t\sigma_3} +\int_{-t}^t B(t,s)
\ee^{-4\ii k^3s\sigma_3}\dd s\notag\\
&\quad
+\ii k\int_{-t}^t C(t,s)\ee^{-4\ii k^3s\sigma_3}\dd s 
+k^2\int_{-t}^t D(t,s)\ee^{-4\ii k^3s\sigma_3}\dd s,
\end{align}
which will be used below. The proof of this triangular representation
can be done by the same way as in \cite{BMK}.
In the present case the matrix-valued real functions
$A(x,y,t)$, $B(t,s)$, $C(t,s)$ and $D(t,s)$ are $C^{\infty}$
and bounded in $x,y,t,s$.

The triangular integral representations \textup{(\ref{PHI})-(\ref{hphi})}
yield the following properties of the solution $\Phi(x,t,k)$:

\subsubsection*{Properties of the 
                second basic solution}
\begin{enumerate}
\item  
$\Phi (x,t,k)$ is a solution of the $x$- and $t$-equations.
\item  
$\Phi (x,t,k)=\Lambda\bar{\Phi}(x,t,\bar{k})\Lambda^{-1}$ for any
$k\in\D{C}$.
\item  
$\Phi (x,t,k)=\bar{\Phi}(x,t,-\bar{k})$ for any 
$k\in\D{C}$.
\item  
$\det \Phi (x,t,k)\equiv 1$ for any $k\in\D{C}$.
\item  
$(x,t,k)\mapsto 
\Phi (x,t,k)\in C^{\infty}(\D{R}_+\times\D{R}_+\times\D{C})$.
\item  
$\Phi(x,t,k)$ is analytic (entire) in $k\in\D{C}$.
\item  
For $k\in\D{C}$, $k\to\infty $,
\[
\Phi (x,t,k)=\Bigl[I+\ord(k^{-1})+\ord\Bigl(\frac{\ee^{2\ii kx\sigma_3}}{k}\Bigr)
+\ord\Bigl(\frac{\ee^{8\ii k^{3}t\sigma_3}}{k}\Bigr)
\Bigr]\ee^{-\ii k(x+4k^2t)\sigma_3}.
\]
\item  For $k\in\Omega_1\cup\Omega_3$, $k\to\infty $,
\[
\ee^{\ii kx+4\ii k^{3}t}\Phi ^{-}(x,t,k)=
\begin{pmatrix}
1 \\
0
\end{pmatrix}
+\ord(k^{-1}).
\]
\end{enumerate}
The last asymptotic relation can be easily proved using large $k$
asymptotics for the functions $\varphi^\mp(x,t,k)$ and $\hat\varphi^-(t,k)$.

\subsection{Third basic solution}

Let $\Sigma=\{k\in\D{C}\mid\Im k^3=0\}$ as above and let $\hat\Psi(t,k)$ be a
solution of the Volterra integral equation
\begin{align}  \notag
\hat\Psi(t,k)&=e^{-4\ii k^3t\sigma_3} -
\int_{t}^{\infty}e^{-4\ii k^3(\tau-t)\sigma_3}\hat Q(0,\tau,k)
\hat\Psi(\tau,k)\dd\tau, \qquad k\in\Sigma,
\end{align}
where $\hat Q(0,t,k)$ is as in (\ref{hatphi}) and $\hat Q(0,t,k)\equiv0$
for $t>T$ if $T<\infty$,
that means the matrix $\hat\Psi(t,k)$ satisfies the $t$-equation with $x=0$
under the asymptotic condition $\hat\Psi(t,k)=e^{-4\ii k^3t\sigma_3}+\R{o}(1)$
as $t\to\infty$.
Again, for $\hat\Psi(t,k)$ the triangular integral representation
\begin{align}                                     \label{hatpsi}
\hat\Psi(t,k)
&=\ee^{-4\ii k^3t\sigma_3}+\int_{t}^{\infty}L(t,s)
\ee^{-4\ii k^3s\sigma_3}\dd s\notag\\
&\quad
+\ii k\int_{t}^{\infty} M(t,s)
\ee^{-4\ii k^3s\sigma_3}\dd s+k^2\int_{t}^{\infty} N(t,s)
\ee^{-4\ii k^3s\sigma_3}\dd s
\end{align}
can be obtained as in \cite{BMK}. Here the matrix-valued real
functions
$L(t,s)$, $M(t,s)$ and  $N(t,s)$ are $C^{\infty}$ and bounded in $t,s$, 
and they vanish for $s>2T-t$ if $T<\infty$. We introduce the matrix
\begin{equation}  \label{Y}
Y(x,t,k)=\varphi(x,t,k)\hat\Psi(t,k),\qquad k\in\Sigma,
\end{equation}
where $\varphi(x,t,k)$ is as in (\ref{phi}). Lemma \ref{lem.2.1} implies
that $Y(x,t,k)$ is a solution of the $x$- and $t$-equations with
\[
\det Y(x,t,k)=1\text{ for }k\in\Sigma.
\]
For $k\notin\Sigma$ the function $\hat\Psi(t,k)$,
hence also $Y(x,t,k)$, is unbounded in $t\in\D{R}_+$.
Since  the integral equation  is of
Volterra type with $\tau\in(t,\infty)$, the first column $Y^-(x,t,k)$ is
analytic in $k\in\Omega_2\cup\Omega_4\cup\Omega_6$ and the second column
$Y^+(x,t,k)$ is analytic in $k\in\Omega_1\cup\Omega_3
\cup\Omega_5$ or $Y(x,t,k)$ is an entire matrix-valued function if $T<\infty$.

The properties of the solution $Y(x,t,k)$ follow from the triangular
integral representations (\ref{phi})  and  (\ref{hatpsi}):

\subsubsection*{Properties of the 
                third basic solution}
\begin{enumerate}
\item  
$Y(x,t,k)$ satisfies the $x$- and $t$-equations.
\item  
$Y(x,t,k)=\Lambda\bar{Y}(x,t,\bar{k})\Lambda^{-1}$ for $k\in\Sigma $.
\item  
$\det Y(x,t,k)=1$ for $k\in\Sigma $.
\item  
$(x,t,k)\mapsto Y(x,t,k)\in C^{\infty}(\D{R}_+\times\D{R}_+\times\Sigma)$.
\item  
$Y^+(x,t,k)$ is analytic in 
$k\in\Omega_1\cup\Omega_3\cup\Omega_5$, $Y^-(x,t,k)$ is analytic in 
$k\in\Omega_2\cup\Omega_4\cup\Omega_6$
or they are entire if $T<\infty$.
\item  
$\bar Y^{-}(x,t,-\bar k)= Y^{-}(x,t,{k})$,\;
$k\in\Omega_2\cup \Omega_4\cup\Omega_6$.
\item  
$\bar Y^{+}(x,t,-\bar k)=Y^{+}(x,t,{k})$,\;
$k\in\Omega_1\cup \Omega_3\cup\Omega_5$.
\item  
For $k\to\infty $,\; $k\in\Omega_2$,
\[
\ee^{\ii kx+4\ii k^{3}t}Y^{-}(x,t,k)=
\begin{pmatrix}
1\\
0
\end{pmatrix}
+\ord(k^{-1}).
\]
\end{enumerate}

\section{Analysis of the direct 
                   scattering problem}
\setcounter{equation}{0}

The basic solutions we have introduced are clearly linearly dependent:
\begin{align}  \label{T}
\Phi(x,t, k)&= \Psi(x,t,k)S(k),  \notag \\
Y(x,t,k)&= \Phi(x,t,k)P(k), \\
Y(x,t,k)&= \Psi(x,t,k)R(k).  \notag
\end{align}
The matrices $S(k)$, $P(k)$ and $R(k)$ depend neither on $x$ nor on $t$
because by virtue of the $x$-equation they do not depend on $x$, and by
virtue of the $t$-equation they do not depend on $t$. Hence:
\begin{alignat}{2}  \label{T1}
S(k)&=\Psi^{-1}(0,0,k), & \quad& k\in\D{R};  \notag \\
P(k)&=Y(0,0,k), & &k\in\Sigma; \\
R(k)&=S(k)P(k), & &k\in\D{R}.  \notag
\end{alignat}
Let us study the properties of these ``scattering'' (transition) matrices.

\subsubsection*{Properties of the 
              scattering matrix $S(k)$}
They follow from the scattering problem for the $x$-equation
with $t=0$. Indeed, consider the problem on the whole $x$-line by putting
\[
q(x,0)=\hat u(x)=
\begin{cases}
0 & \text{ for } x\in(-\infty,0) \\
u(x) & \text{ for } x\in [0,\infty).
\end{cases}
\]
Let $\tilde \Psi(x,k)$ be the Jost solution \cite{FT}
normalized by
\[
\tilde\Psi(x,k)=\ee^{-\ii kx\sigma_3}\text{ for }x<0
\]
and let $\tilde T(k)$ be the transition matrix for that case, i.e.
\[
\tilde\Psi(x,k)= \Psi(x,0,k)\tilde T(k).
\]
Putting $x=0$ we find $S(k)\equiv\tilde T(k)$. Hence the ``scattering''
matrix $S(k)$ has all properties of the transition matrix $\tilde T(k)$ 
\cite{FT}:
\begin{itemize}
\item  
$S(k)=\Lambda\bar{S}(k)\Lambda^{-1}$ \text{for } $k\in\D{R}$.
\item  
$\det S(k)\equiv 1$ \text{for } $k\in\D{R}$.
\item  
$S(k)\in C^{\infty}(\D{R})$.
\end{itemize}
For the half-line case there are additional properties:
\begin{itemize}
\item  
$S(k)=
\begin{pmatrix}
s_2^{+}(k) & -s_1^{+}(k) \\
-s_2^{-}(k) & s_1^{-}(k)
\end{pmatrix}
$ where $s_{j}^{+}(k)=\Psi _{j}^{+}(0,0,k)$;
\item  
$\begin{pmatrix}s_2^{+}(k)&-s_1^{+}(k)\end{pmatrix}$ is analytic in $k\in\D{C}_+$
and $s_j^+(k)=\overline{ s_j^+(-\bar k)}$;
\item  
$\begin{pmatrix}-s_2^{-}(k)&s_1^{-}(k)\end{pmatrix}$ is 
analytic\footnote{For an arbitrary function $\hat{u}(x)$, $x\in\D{R}$ these analytic
properties do not hold, $s_2^{+}(k)$ and $s_1^{-}(k)$ are only analytic
in $k\in\D{C}_+$ and $k\in\D{C}_-$ respectively. In our case
$\hat{u}(x)\equiv 0$ for $x<0$ and therefore $s_1^+(k)$ and $s_2^-(k)$ are
also analytic in $k\in\D{C}_+$ and $k\in\D{C}_-$ respectively.}
in $k\in\D{C}_-$
and $s^-_j(k)=\overline{s^-_j(-\bar k)}$;
\item  
If $k\in\D{C}_+$ and $k\to\infty$,
\begin{align}         \label{s12}
s_2^{+}(k) =1+\ord(k^{-1}), \qquad s_1^{+}(k) =\ord(k^{-1}).
\end{align}
\end{itemize}
Let us prove some integral representations for $s_1^+(k)$ and $s_2^+(k)$.
We use the limit formulas
\begin{align*}
s_1^+(k)=\lim_{x\to-\infty}\ee^{-\ii kx}\Psi_1^+(x,0,k), \quad
s_2^+(k)=\lim_{x\to-\infty}\ee^{-\ii kx}\Psi_2^+(x,0,k),
\end{align*}
which follow from the definition of the matrix $S(k)$. If one puts
\[
\chi_j(x,k)=\ee^{-\ii kx}\Psi_j^+(x,0,k)
\]
the $x$-equation yields
\begin{alignat*}
{2} &\chi_1'+2\ii k\chi_1=\hat u(x)\chi_2
&\qquad&\chi_1(x,k)\to 0\text{ as } x\to+\infty \\
&\chi_2'=-{\hat u(x)}\chi_1&&\chi_2(x,k)\to 1\text{ as }
x\to+\infty.
\end{alignat*}
Then, by integration:
\begin{align*}
\chi_1(x,k)&=-\int_x^{\infty}\ee^{2\ii k(y-x})\hat
u(y)\chi_2(y,k)\dd y \\
\chi_2(x,k)&=1+\int_x^{\infty}{\hat u(y)}\chi_1(y,k)\dd y,
\end{align*}
therefore
\begin{align}             \label{S12}
s_1^+(k)&=-\int^{\infty}_0 u(x)\ee^{2\ii kx}\dd x
-\int^{\infty}_0 u(x)\ee^{2\ii kx}\dd x
\int^{\infty}_0K_1(x,x+y)\ee^{\ii ky}\dd y \\
s_2^+(k)&=1-\int^{\infty}_0\ee^{\ii kx}\dd x
\int^{\infty}_0{ u(y)} K_2(y,y+x)\dd y,
\end{align}
where $K_1(x,y)$ and $K_2(x,y)$ are entries of the kernel of
triangular integral transformation (\ref{PSI}). The last two formulas allow
to find the large-$k$ asymptotic expansions at any order of $s_1^+(k)$ and
$s_2^+(k)$ and to obtain (\ref{s12}), in particular,
which is exact (precise) if $u(0)\neq 0$.

The matrix $S(k)=\Psi^{-1}(0,0,k)$ is determined by 
$u(x)\in\C{S}(\D{R}_+)$. The entries of this matrix are
not independent and
can be recovered from one known function. 
Let 
\[
s(k)\equiv\frac{s_1^+(k)}{s_2^+(k)}
\]
be given and let 
\[
\Sigma_{\R{d}}^{\text{ic}}=\{k_1,\dots,k_n\in\D{C}_+\mid s_2^+(k_j)=0\},
\]
be the set of zeros of the analytic function $s_2^+(k)$,
which is finite because 
$s_2^+(k)\to 1$ as $k\to\infty$ and we have supposed that $s_2^+(k)\neq 0$
for any $k\in\D{R}$.
Since $\det S(k)\equiv 1$, then $|s_2^+(k)|^2-\lambda|s_1^+(k)|^2\equiv 1$ for
any $k\in\D{R}$. This identity yields the well-known formula:
\begin{equation}  \label{s2k}
s_2^+(k)=\left(\prod_{k_j\in{\D{C}_+}}\frac{k-k_j}{k-\bar k_j}\right)
^{\frac{1-\lambda}{2}} 
\exp\left\{\frac{\ii}{2\pi}\int_{-\infty}^{\infty}
\frac{\log[1-\mu|s(\mu)|^2]
\dd\mu}{\mu-k}\right\},
\end{equation}
The remaining entries of $S(k)$ are also recovered:
\[
s_1^+(k)=s(k)s_2^+(k), 
\quad 
s_2^-(k)=\lambda\bar s_1^+(\bar k), 
\quad
s_1^-(k)=\bar{s}_2^+(\bar k).
\]
So, if $\lambda=1$, the function $s_2^+(k)$ has no zeros at all and the set
$\Sigma_{\R{d}}^{\text{ic}}$ is  empty. It follows from the
self-adjointness of the $x$-equation (\ref{xeq}) and the obvious inequality:
$|s_2^+(k)|\geq 1$ for $k\in\D{R}$.
If $\lambda=-1$ then $s_2^+(k)$
may vanish at some points $k_j\in\D{C}_+$. Since
$u(x)$ is real-valued, $\Sigma_{\R{d}}^{\text{ic}}$ is symmetric 
with respect to the imaginary axis. 
We can enumerate the $k_j$'s in
such a way that
\begin{itemize}
\item
$k_j=\ii\kappa_j$,
$\kappa_j>0$ for $j=1,\dots,n_1\leq n$ with $n=n_1+2n_2$,
\item
$k_{n_1+l}=-\bar k_{n_1+n_2+l}$ for $l=1,\dots,n_2$.
\end{itemize}
Moreover, these zeros
can be multiple and there can exist limit points on the real line
$\D{R}$ \cite{FT}. 

To avoid this difficulties we shall consider a
subset $\C{S}_0(\D{R}_+)$ of functions 
$u(x)\in\C{S}(\D{R}_+)$ for which $s_2^+(k)$ has a finite number of
zeros $k_1,\dots,k_n$ in $\D{C}_+$, all of multiplicity $1$, i.e.\ $\dot
s_2^+(k_j)\neq 0$, and $s_2^+(k)\neq 0$ for every $k\in\D{R}$.

Let us briefly discuss the discrete spectrum of the $x$-problem, which
may appear when $\lambda=-1$. The main
relation of the $x$-scattering problem is
\begin{equation}  \label{xscat}
\frac{1}{s_2^+(k)}\Phi^-(x,t,k)=\Psi^-(x,t,k)+r(k) \Psi^+(x,t,k) \text{ for }
k\in\D{R},
\end{equation}
where
\begin{equation}  \label{r}
r(k)=-\frac{s_2^-(k)}{s_2^+(k)}.
\end{equation}
$F(x,t,k)=\Phi^-(x,t,k)/s_2^+(k)$ is analytic in $k\in\D{C}_+$
except for 
$\Sigma_{\R{d}}^{\text{ic}}=\{k_1,\dots,k_n\}$,
where it has poles. We have
\[
s_2^+(k_j)=\det\begin{bmatrix}\Phi^-(x,t,k_j)&\Psi^+(x,t,k_j)\end{bmatrix}=0,
\]
then $\Phi^-(x,t,k_j)=\gamma_j^1\Psi^+(x,t,k_j)$. Hence,
\[
\res_{k=k_j}F(x,t,k)=c_j^1\Psi^+(x,t,k_j)
\]
with
\[
c_j^1=\frac{\gamma_j^1}{\dot s_2^+(k_j)}\quad
\text{and}\quad \gamma_j^1=\frac{1}{s_1^+(k_j)},\quad j=1,\dots,n.
\]
The dot denotes differentiation with respect to $k$. Note that $s_1^+(k_j)
\neq0$ because otherwise we come to a contradiction: $\Psi_+(x,t,k_j)\equiv 0
$ since $\Psi_1^+(0,0,k_j)=s_1^+(k_j)=0$ and $\Psi_2^+(0,0,k_j)=s_2^+(k_j)=0$.
 We also assume all zeros are simple, i.e.\ 
$\dot s_2^+(k_j)\neq 0$.

Using asymptotics of the function $\Phi^-(x,t,k)$ at $k=\infty$,
for $k\in\Omega_1\cup\Omega_3$, we find
\begin{equation}  \label{Fas}
F(x,t,k)=
\left[
\begin{pmatrix}
1 \\
0
\end{pmatrix}
+ \ord(|k|^{-1})\right]\ee^{-\ii kx-4\ii k^3t}
\text{ for }|k|\to\infty,\;  k\in\Omega_1\cup\Omega_3,
\end{equation}
that will be used below. So, we come to the following (cf.\ conditions
A and B):

\subsubsection*{Properties of 
$\BS{r(k)}$, $\BS{t(k)}$ and $\BS{k_j}$}
\begin{itemize}
\item  
The reflection coefficient $r(k)$ belongs to 
$C^{\infty}(\D{R})$, $r(-k)=\bar r(k)$ and $r(k)=\ord(k^{-1}) $ as
$k\to\infty$. It is the ratio of two functions $-s_2^-(k)$ 
and $s_2^+(k)$ analytic in $k\in\D{C}_-$ and $k\in\D{C}_+$
respectively, and $|r(k)|<1$ if $\lambda=1$.
\item  
The transition coefficient $t(k)=[s_2^+(k)]^{-1}$ is represented
through formula (\ref{s2k}), where $|s(\mu)|=|r(\mu)|$.
\item  
If $\lambda=-1$, $k_i\neq k_j$ for $i\neq j$ and $\Im k_j>0$, 
$j=1,\dots,n$, with:\\  
$k_j=\ii\kappa_j$, $1\leq j\leq n_1\leq n=n_1+2n_2$, and
$k_{n_1+l}=-\bar k_{n_1+n_2+l}$, $1\leq l\leq n_2$.
\end{itemize}

These properties follow from the $x$-scattering problem on the whole line.
We take into account  that
for the half-line case the function $s_2^-(k)=-\bar s_1^+(\bar k)$ is
analytic in $k\in\D{C}_-$ and the constants $\gamma_j^1$ are not
independent parameters (that takes place for the whole line), and they are
evaluated by means of the
function $s_1^+(k)$ at $k_j$: $\gamma_j^1=1/s_1^+(k_j)$.

\subsubsection*{Properties of the 
              scattering matrix $\BS{P(k)}$}
They follow from the
defining relation, i.e.\ from (\ref{T1}), (\ref{Y}), (\ref{hatpsi}):
\[
P(k)=I +\int_0^{\infty} L(t,s)
\ee^{-4\ii k^3s\sigma_3}\dd s+k\int_0^{\infty} M(t,s)
\ee^{-4\ii k^3s\sigma_3}\dd s +k^2\int_0^{\infty} N(t,s)
\ee^{-4\ii k^3s\sigma_3}\dd s.
\]
It is easy to find the following properties:
\begin{itemize}
\item  
$P(k)=\Lambda\bar{P}(\bar{k})\Lambda^{-1}$ for $k\in\Sigma $.
\item  
$\det P(k)\equiv 1$ for $k\in\Sigma$.
\item  
$P(k)$ is $C^{\infty}$ in $k\in\Sigma$.
\item  
If $T<\infty$ the matrix-valued function $P(k)$ is entire in $k\in\D{C}$.
\item  
If $T=\infty$ the vector-function $P^+(k)$ is analytic 
in $k\in\Omega_1\cup \Omega_3\cup\Omega_5$, and $P^-(k)$ is analytic
in $k\in\Omega_2\cup\Omega_4\cup\Omega_6$.
\item 
$P^{\pm }(k)=\overline{P^{\pm }(-\bar{k})}$.
\item  
$P(k)=\sigma_0+\ord(k^{-1})$,\; $k\in\Sigma $, $k\to\infty$.
\end{itemize}

\subsubsection*{Properties of the 
              scattering matrix $\BS{R(k)}$}
Below we need to study properties of the ``scattering'' matrix $R(k)$
introduced by equations (\ref{T})-(\ref{T1}). The matrix $R(k)$ has the
following form
\begin{align*}
&R(k)=
\begin{pmatrix}
r_1^-(k) & r_1^+(k) \\
r_2^-(k) & r_2^+(k)
\end{pmatrix},
\quad
r_1^+(k)=-\bar r_2^-(\bar k), 
\quad 
r_2^+(k)=\bar r_1^-(\bar k)\
\text{ for } k\in\Sigma
\end{align*}
with
\begin{equation}  \label{r1}
r_1^-(k)= p_1^-(k)s_2^+(k)- p_2^-(k)s_1^+(k)
\end{equation}
analytic in $k\in\D{C}_+$ (if $T<\infty$) and
$k\in\Omega_2$ (if $T=\infty$). Hence $r_2^+(k)$ is analytic in
$k\in\D{C}_-$ (if $T<\infty$) and $k\in\Omega_5$ (if $T=\infty$).
Furthermore
\begin{equation}  \label{r2}
r_2^-(k)= p_2^-(k)s_1^-(k)- p_1^-(k)s_2^-(k)
\end{equation}
is analytic in $k\in\D{C}_-$ (for $T<\infty$) and
$k\in\Omega_4\cup\Omega_6$ (for $T=\infty$),
hence $r_1^+(k)$ is analytic in $k\in\D{C}_+$ (for $T<\infty$) and
$k\in\Omega_1\cup\Omega_3$ (for $T=\infty$).
In the domains of analyticity we have the following symmetry properties~:
\[
r^\pm_1(k)=\overline{r^\pm_1(-\bar k)}\quad \text{and}\quad
r^\pm_2(k)=\overline{r^\pm_2(-\bar k)}.
\]
From (\ref{T}) we derive
\[
Y^+(x,t,k)=r_1^+(k)\Psi^-(x,t,k) + r_2^+(k)\Psi^+(x,t,k)
\]
with
\begin{align*}
&r_1^+(k)=\det\begin{bmatrix}Y^+(x,t,k)&\Psi^+(x,t,k)\end{bmatrix},
&r_2^+(k)=\det\begin{bmatrix}\Psi^-(x,t,k)& Y^+(x,t,k)\end{bmatrix}.
\end{align*}
Let us put $x=0$, and $k=k_1+\ii k_2\in\Omega_1\cup\Omega_3$. Using (\ref
{PSI}), (\ref{Y}), (\ref{hatpsi}) for large $t$ we obtain
\[
|r_1^+(k)|\leq C_1(k)\exp\bigl[8t\bigl(K-(3k^2_1-k^2_2)k_2\bigr)\bigr]
\]
where $C_1(k)$ is independent of $t$, and
\[
K=\max_{1\leq j\leq n}\bigl[(3\Re^2 k_j-\Im^2 k_j)\Im k_j\bigr],
\]
where $k_j\in\Omega_1\cup\Omega_3$ is an eigenvalue of the $x$-scattering
problem. Taking into account the analyticity of the function $r_1^+(k)$ for
$k\in\Omega_1\cup\Omega_3$, choosing a large enough $k$ and putting $t\to\infty$
(if $T=\infty$) we find $r_1^+(k)\equiv 0$ for any
$k\in\Omega_1\cup\Omega_3$, hence $r_2^-(k) \equiv 0$ for any
$k\in\Omega_4\cup\Omega_6$. So, we come to the
main property of the compatible scattering problem for $x$- and 
$t$-equations. 
\begin{itemize}
\item
If $T=\infty$ the ``scattering'' matrix $R(k)$ is diagonal:
\[
R(k)=
\begin{pmatrix}
\rho_-(k) & 0 \\
0 & \rho_+(k)
\end{pmatrix}
\text{ for } k\in\D{R}
\]
with
\begin{equation}  \label{rhopm}
\rho_+(k)=\frac{p_2^+(k)}{s_2^+(k)}=\frac{p_1^+(k)} {s_1^+(k)}, \qquad
\rho_-(k)=\frac{p_1^-(k)}{s_1^-(k)}=\frac{p_2^-(k)} {s_2^-(k)}.
\end{equation}
\item
We also have the important relation
\begin{equation}  \label{p/s}
\frac{p_1^+(k)}{p_2^+(k)}\equiv\frac{s_1^+(k)}{s_2^+(k)}
\end{equation}
that says: the function $p(k):=p_1^+(k)/p_2^+(k)$ being meromorphic in
the domain $\Omega_1\cup\Omega_3$ has an analytic continuation into $\D{C}_+$
up to the meromorphic function $s(k)=s_1^+(k)/s_2^+(k)$.
\item
If $T<\infty$, instead of (\ref{p/s}) we have
the  so-called \emph{global relation} \cite{Fo}~:
\[
s_2^+(k)p_1^+(k,T)-s_1^+(k)p_2^+(k, T)=r_1^+(k, T),
\]
where we write down the dependence on $T$ to emphasize
that functions $p_j^+(k)$ ($j=1,2$) and $r_1^+(k)$
really depend on $T$.
\end{itemize}
The asymptotic behavior of $S(k)$ and $P(k)$ yields the following asymptotic
expansions~:
\[
r_1^-(k)=1+\frac{\rho_1}{k}+\ldots,\quad
r_2^-(k)=\frac{\omega_1}{k}+\frac{\omega_2}{k^2}+\ldots \quad
\text{for }k\to\pm\infty.
\]
If $T=\infty$ then
$r_2^-(k)\equiv r_1^+(k)\equiv0$.
Since $\det R(k)\equiv1$, then $\rho_-(k)\rho_+(k)=|\rho_+(k)|^2\equiv1$.
Hence, $\rho_\pm(k)$ can be written in the form:
\begin{equation}  \label{einu}
\rho_\pm(k)=\ee^{\pm\ii\nu(k)} \qquad k\in\D{R}
\end{equation}
with a real function $\nu(k)$ for $k\in\D{R}$. The function $\nu(k)$
has an analytic continuation to the domain 
$\Omega_1\cup\Omega_3\cup\Omega_4\cup\Omega_6$ which satisfies:
\begin{itemize}
\item
$\nu(k)=\bar\nu(\bar k)$,\; $\nu(k)=-\overline{\nu(-\bar k)}$, and
$\nu(k)\to 0$ as $k\to\infty$,
\end{itemize}
in view of the
asymptotics of the functions $s_2^+(k)$ and $p_2^+(k)$.
Indeed, in view of (\ref{rhopm}),
\begin{equation}
p_j^+(k)=\rho_+(k)\  s_j^+(k)\text{ for }j=1,2,  \label{p=rhos}
\end{equation}
$\rho_+(k)$ must have poles at the points where $s_1^+(k)$ and
$s_2^+(k)$ vanish. On the other hand, $s_1^+(k)$ and $s_2^+(k)$
must simultaneously vanish at poles in view of the analyticity of the
functions $p_j^+(k)$ for $k\in\Omega_1\cup\Omega_3$.
Hence $\Psi^+(x,t,k)$ vanish identically if $k$ is a
pole, which is impossible. So $\rho_+(k)$ is analytic (without
singularities) in $k\in\Omega_1\cup\Omega_3$. Hence
the functions $p_j^+(k)$ and $s_j^+(k)$ have a common set of zeros,
possibly empty, in $\Omega_1\cup\Omega_3$. The other statements
about $\nu(k)$ are obvious. 

So, for $r_1^-(k)$ which
is analytic in $k\in\Omega_2$ we obtain:
\begin{equation}
r_1^-(k)=\frac{s_2^{+}(k)}{p_2^{+}(k)}=
\ee^{-\ii\nu (k)}
\text{ for }k\in {\partial\Omega_2}, \label{rho|i}
\end{equation}
The last formula follows from relations:
\begin{align*}
& r_1^-(k)=p_1^{-}(k)s_2^{+}(k)-p_2^{-}(k)s_1^{+}(k)
& p_1^{+}(k)s_2^{+}(k)-p_2^{+}(k)s_1^{+}(k)=0.
\end{align*}
Hence the function $r_1^-(k)$ has an analytic continuation to the domain
$\Omega_1\cup\Omega_3$, where it coincides with the function
$1/\rho_{+}(k)$. Therefore the function $r_1^-(k)$ does not vanish for
$k\in\partial\Omega_2$, and
its zeros are some points $z_j\in\Omega_2$.
Let $\Sigma_{\R{d}}^{\text{bc}}$ be the set of zeros of the function $r_1^-(k)$. As above we
assume that the number of zeros is finite:
\[
\Sigma_{\R{d}}^{\text{bc}}=\{z_1,\dots,z_m\in\Omega_2\mid r_1^-(z_j)=0\}.
\]
We also assume all zeros are simple, i.e.\ $\dot{r}_1^-(z_j)\neq 0$.
We have $\Sigma_d^{bc}=\varnothing$ if $\lambda=1$. Let
\[
\rho(k):=\frac{r_2^-(k)}{r_1^-(k)}.
\]
The functions $\rho(k)$ and $r_1^-(k)$ are
dependent. They satisfy the determinant relation
\begin{equation}    \label{ro1}
1-\lambda|\rho(k)|^{2}=\frac{1}{|r_1^-(k)|^{2}},\quad k\in\D{R},
 \quad 
\end{equation}
and $\rho(k)\equiv 0$ if $T=\infty$.
They have the following properties:
\begin{itemize}
\item  
$\rho(k)$, $r_1^-(k)\in C^{\infty}(\D{R})$, and
$\rho(-k)=\bar\rho(k)$,\; $r_1^-(k)=\bar r_1^-(-k)$,
\item 
If $T=\infty$, $\rho(k)\equiv 0$ for $k\in\D{R}$, and
$r_1^-(k)=\ee^{-\ii\nu(k)}$ where $\nu(k)$ is described above.
\item 
If $T<\infty$, then
\[
r_1^-(k)=
\biggl(\prod_{z_j\in\D{C}_+}\dfrac{k-z_j}{k-\bar{z}_j}
\biggr)^{\!\!\frac{1-\lambda}{2}}\!
\exp\Bigl[\frac{\ii}{2\pi}
     \int_{-\infty}^{\infty}
     \dfrac{\log (1-\lambda|\rho (s)|^{2})\dd s}{s-k}\Bigr],\quad k\in\D{C}_+.
\]
\item  
The function
\[
\rho(k)-r(k)=\frac{p_2^{-}(k)}{r_1^-(k)s_2^{+}(k)}
\]
has an analytic continuation to $\D{C}_+$ for $T<\infty$.\\
For $T=\infty$ the r.h.s.\ is analytic only in $\Omega_2$.
\end{itemize}
The last item follows
from equations (\ref{r1}) and (\ref{r2}) which yield
\begin{align*}
p_1^-(k)&=r_2^-(k)s_1^+(k)+r_1^-(k)s_1^-(k)=
r_1^-(k)[1/s_2^+(k)+s_1^+(k)(\rho(k)-r(k))], \\[1mm]
p_2^-(k)&=r_2^-(k)s_2^+(k)+r_1^-(k)s_2^-(k)=
r_1^-(k)s_2^+(k)[\rho(k)-r(k)] \text{ for } k\in\D{R}.
\end{align*}
For $T<\infty$ the difference $\rho(k)-r(k)$ has an analytic continuation to
$\D{C}_+$  because the l.h.s.\ are analytic in $k\in\D{C}_+$. 
Hence, the r.h.s.\ must have
analytic continuations to $\D{C}_+$.

The second main relation of the compatible scattering problem is:
\begin{align}
G(x,t,k)
&=\frac{1}{r_1^-(k)}Y^{-}(x,t,k)    \label{xscat2}\\
&=
\begin{cases}
\Psi ^{-}(x,t,k)+\rho(k)\Psi^{+}(x,t,k)\
&\text{for }k\in\D{R}\text{ if }T<\infty\\
\Psi^-(x,t,k)
&\text{for }k\in\D{R}\text{ if }T=\infty.
\end{cases}\notag
\end{align}
The function $G(x,t,k)$ is analytic in $k\in\Omega_2$, for
$k\neq z_j$ and the $z_j$'s are poles of $G$. If $r_1^-(z_j)=0$ then
$Y^{-}(x,t,z_j)$ and $\Psi^+(x,t,z_j)$ are linearly dependent:
\[
Y^{-}(x,t,z_j)=\gamma_j^2\Psi^+(x,t,z_j),\quad j=1,\dots,m,
\]
hence
\[
\res_{k=z_j}G(x,t,k)=c_j^{2}\Psi^+(x,t,z_j),\qquad c_j^{2}=
\frac{\gamma _j^{2}}{\dot{r}_1^-(z_j)}
\]
(the dot denotes differentiation with respect to $k$) with
\[
\gamma_j^{2}=\frac{p_1^{-}(z_j)}{s_1^{+}(z_j)}=
\frac{p_2^{-}(z_j)}{s_2^{+}(z_j)}.
\]
Using asymptotics of the function $Y^{-}(x,t,k)$ in the neighborhood of
$k=\infty$ for $k\in\Omega_2$, we find
\begin{align}
& G(x,t,k)=
\left[
\begin{pmatrix}
1 \\
0
\end{pmatrix}
+\ord(|k|^{-1})\right] \ee^{-\ii kx-4\ii k^{3}t}
\label{Gas} \\
& \text{for }|k|\to\infty ,\ k\in\Omega_2.  \notag
\end{align}
This asymptotic formula will be used in the next section.

\section{The main integral equations}
\setcounter{equation}{0}

The main relations of the compatible scattering problem follow from 
(\ref{T}), (\ref{T1}) and (\ref{xscat}), (\ref{xscat2}):
\begin{align}  \label{F}
F(x,t,k)&=\Psi^-(x,t,k)+r(k) \Psi^+(x,t,k) \text{ for }
k\in\D{R}, \\\label{G}
G(x,t,k)&=\Psi^-(x,t,k)+ \rho(k)\Psi^+(x,t,k) \text{ for } k\in\D{R}.
\end{align}
These relations give:
\begin{equation}  \label{c0}
G(x,t,k) -F(x,t,k)=c(k)\Psi^+(x,t,k),
\end{equation}
where $c(k)$ can be written as follows
\begin{equation}  \label{c}
c(k)=\rho(k)-r(k)=
\frac{p_2^-(k)}{s_2^+(k)\ r_1^-(k)}\quad 
\text{for }
\begin{cases}
k\in\D{C}_+&\text{if }T<\infty,\\
k\in\Omega_2&\text{if }T=\infty.
\end{cases}
\end{equation}

\subsubsection*{Properties of 
                 $\BS{c(k)}$}
Indeed, $c(k)$ is  meromorphic in
$\D{C}_+$ if $T<\infty$, and in $\Omega_2$ if $T=\infty$, and
$c(k)=-\overline{c(-\bar k)}$,
because $p_2^-(k)$, $s_2^+(k)$ and $r_1^-(k)$ are analytic in
$k\in\D{C}_+$ if $T<\infty$, and in $k\in\Omega_2$ if $T=\infty$.
Hence relation (\ref{c0}) is true for all $k\in\overline\Omega_2$.
The function $c(k)$ has poles at the points $z_j$, where $s_2^+(z_j)=r_1^-(z_j)=0$.
Since the zeros of $s_2^+(k)$ and $r_1^-(k)$ are simple and in finite number,
all poles of $c(k)$ are simple and also in finite number. Indeed, we only
have to check the case
$s_2^+(z_0)=r_1^-(z_0)=0$. Due to (\ref{r1}) we also find $p_2^-(z_0)=0$.

We have the following relation on $\partial\Omega_2$:
\begin{equation}  \label{iR+}
\frac{Y^-(x,t,k-0)}{r_1^-(k-0)} - \frac{\Phi^-(x,t,k+0)}{s_2^+(k+0)}
=c(k)\Psi^+(x,t,k) \text{ for } k\in\partial\Omega_2.
\end{equation}
To deduce the integral equations of the inverse scattering problem let us
put
\begin{align*}
h^{-}(x,t,k)& =G(x,t,k)-
\begin{pmatrix}
1 \\
0
\end{pmatrix}
\ee^{-\ii kx-4\ii k^{3}t}\text{ for }k\in\partial\Omega_2\\
h^{+}(x,t,k)& =F(x,t,k)-
\begin{pmatrix}
1 \\
0
\end{pmatrix}
\ee^{-\ii kx-4\ii k^{3}t}\text{ for }k\in\D{R}.
\end{align*}
Let us consider the integral
\[
J(x,y,t)=\frac{1}{2\pi}\int_{\partial\Omega_2} h^-(x,t,k)
\ee^{\ii ky+4\ii k^3t}\dd k +\frac{1}{2\pi}
\int^{\infty}_{-\infty} h^+(x,t,k)\ee^{\ii ky+4\ii k^3t}\dd k.
\]
Using equations (\ref{F}), (\ref{G}), (\ref{c}), (\ref{PSI}) we find
\begin{align*}
&J(x,y,t)-\frac{1}{2\pi}\int_{\partial\Omega_2} h^-(x,t,k)
\ee^{\ii ky+4\ii k^3t}\dd k=\\
&\qquad\qquad
=
\begin{pmatrix}
K_1 \\
K_2
\end{pmatrix}
(x,y,t)+
\begin{pmatrix}
0 \\
1
\end{pmatrix}
F_s(x+y,t) 
+\int_x^{\infty}
\begin{pmatrix}
\lambda K_2 \\
K_1
\end{pmatrix}
(x,z,t) F_s(z+y,t)\dd z,
\end{align*}
where
\begin{align*}
F_s(x,t)&=\frac{1}{2\pi}\int_{-\infty}^{\infty} r(k)
\ee^{\ii k(x+y)+8\ii k^3t}\dd k.
\end{align*}
On the other hand, using estimates (\ref{Fas}) and (\ref{Gas}) of $F(x,t,z)$
and $G(x,t,z)$ for large $k$, taking into account 
(\ref{c0}) (\ref{c}), (\ref{iR+}) and
applying the Jordan lemma, we find
\begin{align*}
J(x,y,t)
&=\frac{1-\lambda}{2}
\biggl(
\ii\sum_{\begin{subarray}{l}
          k_j\in\Omega_1\cup\Omega_3\\
          s_2^+(k_j)=0
          \end{subarray}}
          \res_{k=k_j}\bigl[h^+(x,t,k)\ee^{\ii ky+4\ii k^3t}\bigr]\\[-2mm]
&\qquad\qquad\qquad\qquad
+\ii\sum_{\begin{subarray}{l}
           z_j\in\Omega_2\\
           r_1^-(z_j)=0
          \end{subarray}}
          \res_{k=z_j}\bigl[h^-(x,t,k)\ee^{\ii ky+4\ii k^3t}\bigr]\biggr)\\
&\quad
-\frac{1}{2\pi}\int_{\partial\Omega_2}[G(x,t,k)-F(x,t,k)]\ee^{\ii ky+4\ii k^3t}
\dd k\\
&=
\frac{1-\lambda}{2}
\biggl(-\sum_{k_j\in\Omega_1}
         m_j^{1}\ee^{\ii k_jy+4\ii k_j^3t}\Psi^+(x,t,k_j)
       -\sum_{k_j\in\Omega_3}
         m_j^{3}\ee^{\ii k_jy+4\ii k_j^3t}\Psi^+(x,t,k_j)\biggr)\\
&\;\;\;
-\frac{1-\lambda}{2}
        \sum_{z_j\in\Omega_2}
          m_j^2\ee^{\ii z_jy+4\ii z_j^3t}\Psi^+(x,t,z_j)
 -\frac{1}{2\pi}
         \int_{\partial\Omega_2}c(k)\ee^{\ii ky+4\ii k^3t}\Psi^+(x,t,k)\dd k.
\end{align*}
Finally we have the following integral equations 
of the inverse scattering:
\begin{align}                                              \label{BIE1}
&K_1(x,y,t)+\lambda\int_x^{\infty} K_2(x,z,t)H(z+y,t)\dd z=0 \text{ for }
0\leq x<y<\infty, \\ \label{BIE2}
&K_2(x,y,t)+ H(x+y,t)+\int_x^{\infty}K_1(x,z,t) H(z+y,t)\dd z=0
\end{align}
with the kernel
\begin{align}	                                          \label{H}
H(x,t)
&=\frac{1-\lambda}{2}
\biggl(\sum_{k_j\in\Omega_1}
m_j^1\ee^{\ii k_jx+8\ii k_j^3t}
+\sum_{z_j\in\Omega_2}
m_j^2\ee^{\ii z_jx+8\ii z_j^3t}+
\sum_{k_j\in\Omega_3}
m_j^3\ee^{\ii k_jx+8\ii k_j^3t}\biggr)
\notag\\
&\quad
+\frac{1}{2\pi}\int_{\partial\Omega_2}
c(k)\ee^{\ii kx+8\ii k^3t}\dd k
+\frac{1}{2\pi}\int_{-\infty}^{\infty}r(k)\ee^{\ii kx+8\ii k^3t}\dd k.
\end{align}
The coefficients $m_j^1$ and $m_j^2$ are given by
\begin{align}  \label{m}
m_j^1&=[\ii s_1^+(k_j)\dot s_2^+(k_j)]^{-1},  \notag \\
m_j^2&=p_1^-(z_j) [\ii s_1^+(z_j)\dot r_1^-(z_j)]^{-1} =p_2^-(z_j) 
[\ii s_2^+(z_j)\dot r_1^-(z_j)]^{-1}=-\ii\res_{k=z_j}c(k).
\end{align}
Using (\ref{BIE2}) for $y=x$ one can prove that
$H(x,t)\in C^{\infty}(\D{R}_+\times\D{R}_+)$ and is rapidly decreasing in $x$, i.e.\
$H(x,t)=\ord(x^{-\infty})$, as $x\to\infty$, since $K_1(x,y,t)$ and $K_2(x,y,t)$ are
in $C^{\infty}(\D{R}_+\times\D{R}_+\times\D{R}_+)$
and also rapidly decreasing as $x+y\to\infty$,
and (\ref{BIE2}) is a Volterra integral equation with respect to
the kernel $H(x,t)$. All terms in (\ref{H}) are clearly $C^{\infty}$ 
and of the Schwartz type except for the last term~:
\[
\frac{1}{2\pi}\int_{-\infty}^{\infty}r(k)\ee^{\ii kx+8\ii k^3t}\dd k
\]
because the reflection coefficient $r(k)$ vanish at infinity as
$\ord(k^{-1})$ that follows from (\ref{s12}), (\ref{S12}).
Thus we arrive to the requirement that this integral must be $C^{\infty}$ 
in $x,t$ (third item in condition A).
It is easy to see that under additional condition (\ref{dn0})
the reflection coefficient $r(k)\in\C{S}(\D{R})$, therefore
the third item can be omitted because the corresponding integral
is $C^{\infty}$.
By assumption the boundary functions
$q(0,t)$, $q_x'(0,t)$, $q_{xx}''(0,t)$ are $C^{\infty}$, then we obtain the following condition:
the kernel
$H(x,t)$ must be $C^{\infty}$ in $t$ when $x=0$. Therefore  the function
\[
\int_{\partial\Omega_2}
c(k)\ee^{8\ii k^3t}\dd k
+\int_{-\infty}^{\infty}r(k)\ee^{8\ii k^3t}\dd k.
\]
should be $C^{\infty}$ with fast decay
as $t\to\infty$.  Thus we arrive to the properties given in the fifth item in
condition C.

For any fixed $t\in\D{R}_+$ the function $H(x,t)$ will be rapidly
decreasing for $x\to\infty$.
Indeed, using the method of steepest descent and integration by parts we see
that $H(x,t)=\ord(x^{-\infty})$ because $c(k)$  and  $r(k)$
are $C^{\infty}$
and, according to their asymptotic behavior, they vanish at infinity as well
as their derivatives of any order.

\begin{rem*}
For $t=0$ the kernel $H(x,t)|_{t=0}$ coincides with the kernel
\[
H_0(x)=\sum_{k_j\in\D{C}_+}\  m_j\ee^{\ii k_jx}+\frac{1}{2\pi}
\int_{-\infty}^{\infty}r(k)\ee^{\ii kx}\dd k,
\]
because in this case ($t=0$) the integral over $\partial\Omega_2$ can be evaluated by
using the residues of the function $c(k)$. After integration we find that
$H(x,0)=H_0(x)$. Then Marchenko integral equations with kernel $H_0(x)$
yield $q(x,0)=u(x)$.
\end{rem*}

Now it is natural to introduce the set
\begin{equation}                                   \label{calR}
\C{R}=\{k_1, k_2, \dots, k_n\in\D{C}_+;
\ z_1, z_2, \dots, z_m\in\Omega_2;
\ r(k),\ k\in\D{R}; \ c(k), \ k\in\Omega_2 \}
\end{equation}
and to call it (see \S1.4) the set of scattering data of the compatible eigenvalue problem for
the pair of differential equations (\ref{xeq}), (\ref{teq}) defined by $q(x,t)$
satisfying (\ref{mkd}).
The kernel $H(x,t)$ of the Marchenko equations is completely defined by the
scattering data $\C{R}$ because the three missing coefficients
$m_j^{1}$, $m_j^{2}$, $m_j^3$ (\ref{m})
can be evaluated from the scattering data.
Conditions A\textendash B\textendash C from the Introduction follow immediately 
from the properties proved above for the scattering data $\C{R}$.

The properties given by conditions A\textendash B\textendash C are
characteristic, i.e.\ they are sufficient to ensure that the system of numbers
$k_1,\dots,k_n$, $z_1,\dots,z_m$
and functions $r(k)$, $c(k)$, are the scattering data of
compatible $x$- and $t$-equations (\ref{xeq}), (\ref{teq}) with $q(x,t)$
satisfying the mKdV equation (\ref{mkd}) with an initial
function $u(x)\in\C{S}(\D{R}_+)$ and boundary values
$v(t), v_1(t), v_2(t)\in C^{\infty}[0,T]$ if $T<\infty$, or
$v(t), v_1(t), v_2(t)\in\C{S}(\D{R}_+)$ if $T=\infty$.

Anyway formula (\ref{q=K2}) and
Marchenko integral equations (\ref{BIE1}), (\ref{BIE2}) represent a solution
of the  mKdV equation if the kernel (\ref{H})
is sufficiently smooth and rapidly decreasing as $x\to\infty$. It
follows from statements in Section 6.

\section{Formulation of the Riemann-Hilbert problem}
\setcounter{equation}{0}

Here we give a formulation of the inverse scattering problem as
a Riemann-Hilbert problem which will be used for
proving that the solution $q(x,t)$ arising from Marchenko equations
satisfies the boundary conditions. We recall that the equality
$q(x,0)=u(x)$ is already proved.

The main scattering relations (\ref{T}) yield the following
Riemann-Hilbert problem. Indeed, from (\ref{T}) we derive~:
\begin{alignat}{2}                                \label{IST1}
&\frac{\Phi^-(x,t,k)}{s_2^+(k)}
         =\Psi^-(x,t,k)+r(k) \Psi^+(x,t,k),             &\quad& k\in\D{R},\notag\\
&\frac{Y^-(x,t,k)}{r_1^-(k)}
         =\Psi^-(x,t,k)+ \rho(k)\Psi^+(x,t,k),          && k\in\D{R},\\
&\frac{\Phi^+(x,t,k)}{s_1^-(k)}
         =\Psi^+(x,t,k)+ \lambda\bar r(k)\Psi^-(x,t,k), && k\in\D{R},\notag\\
&\frac{Y^+(x,t,k)}{r_2^+(k)}
         =\Psi^+(x,t,k)+\lambda\bar\rho(k)\Psi^+(x,t,k),&& k\in\D{R},\\
&\frac{Y^-(x,t,k)}{r_1^-(k)}-\frac{\Phi^-(x,t,k)}{s_2^+(k)}
         = c(k)\Psi^+(x,t,k),                           &&k\in\partial\Omega_2,\\
&\frac{Y^+(x,t,k)}{r_2^+(k)}-\frac{\Phi^+(x,t,k)}{s_1^-(k)}
         =\lambda\bar c(\bar k)\Psi^-(x,t,k),&& k\in\partial\Omega_5.\label{IST2}
\end{alignat}
Let us define the sectionally meromorphic (analytic for $\lambda=1$)
matrix $M(k,x,t)$~:
\[
M(k,x,t)=
\begin{cases}
\begin{pmatrix}
\DS{\frac{\Phi_1^-(x,t,k)\ee^{\ii\theta}}
{s_2^+(k)}}& \Psi_1^+(x,t,k)\ee^{-\ii\theta}\\
\DS{\frac{\Phi_2^-(x,t,k)\ee^{\ii\theta}}{s_2^+(k)}}
&\Psi_2^+(x,t,k)\ee^{-\ii\theta}
\end{pmatrix}&
k\in\Omega_1\cup\Omega_3\\
\begin{pmatrix}
\DS{\frac{Y_1^-(x,t,k)\ee^{\ii\theta}}{r_1^-(k)}}
&\Psi_1^+(x,t,k)\ee^{-\ii\theta}\\
\DS{\frac{Y_2^-(x,t,k)\ee^{\ii\theta}}{r_1^-(k)}}
&\Psi_2^+(x,t,k)\ee^{-\ii\theta}
\end{pmatrix}&
k\in\Omega_2\\
\begin{pmatrix}
\Psi_1^-(x,t,k)\ee^{\ii\theta}&\DS{\frac{Y_1^+(x,t,k)\ee^{-\ii\theta}}{r_2^+(k)}}\\
\Psi_2^-(x,t,k)\ee^{\ii\theta}&\DS{\frac{Y_2^+(x,t,k)\ee^{-\ii\theta}}{r_2^+(k)}}
\end{pmatrix}&
k\in\Omega_5\\
\begin{pmatrix}\Psi_1^-(x,t,k)\ee^{\ii\theta}&
\DS{\frac{\Phi_1^+(x,t,k)\ee^{-\ii\theta}}{s_1^-(k)}}\\
\Psi_2^-(x,t,k)\ee^{\ii\theta}&
\DS{\frac{\Phi_2^+(x,t,k)\ee^{-\ii\theta}}{s_1^-(k)}}
\end{pmatrix}&
k\in\Omega_4\cup\Omega_6.
\end{cases}
\]
Here $\theta=\theta(x,t,k)=k(x+4k^2t)$. Then we have the following
Riemann-Hilbert problem
\begin{equation}             \label{RH}
M_-(k,x,t)=M_+(k,x,t)J(k,x,t), \qquad k\in\Sigma
\end{equation}
on the contour $\Sigma=\{k\mid\Im k^3=0\}$.
The orientation on the contour $\Sigma$ is chosen in such a way that the
sign ``$+$'' (resp.\ ``$-$'') corresponds to the left (resp.\ right) boundary values
of the matrix $M(k,x,t)$ in the domains $\Omega_1$, $\Omega_3$, $\Omega_5$
marked by ``$+$'', and in the domains $\Omega_2$, $\Omega_4$, $\Omega_6$
marked by ``$-$''. The corresponding graph is depicted on Figure 2.
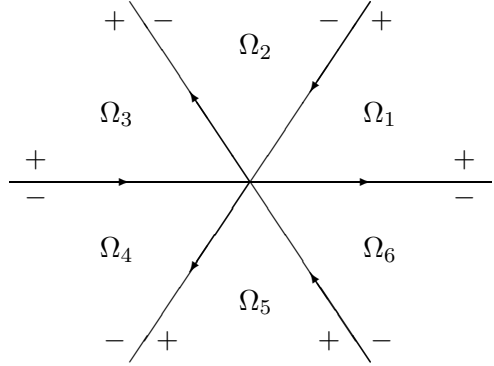
\begin{figure}[ht]
\setlength{\unitlength}{1mm}
\begin{picture}(120,55)
\thinlines
\put(70,30){\line(-2,-3){16}}
\put(70,30){\line(2,3){16}}
\put(70,30){\vector(-2,-3){8}}
\put(86,54){\vector(-2,-3){8}}
\put(70,30){\line(-2,3){16}}
\put(70,30){\line(2,-3){16}}
\put(70,30){\vector(-2,3){8}}
\put(86,6){\vector(-2,3){8}}
\put(38,30){\line(1,0){64}}
\put(70,30){\line(-1,0){10}}
\put(38,30){\vector(1,0){16}}
\put(70,30){\vector(1,0){16}}
\put(85,38){$\Omega_1$}
\put(85,20){$\Omega_6$}
\put(50,38){$\Omega_3$}
\put(68.5,47){$\Omega_2$}
\put(50,20){$\Omega_4$}
\put(68.5,13){$\Omega_5$}
\put(97,32){$+$}
\put(97,27){$-$}
\put(86,50.5){$+$}
\put(79,50.5){$-$}
\put(50.5,50.5){$+$}
\put(57,50.5){$-$}
\put(40,32){$+$}
\put(40,27){$-$}
\put(50.5,8){$-$}
\put(57.5,8){$+$}
\put(86,8){$-$}
\put(79,8){$+$}
\end{picture}
\caption{The oriented contour $\Sigma$.}
\label{fig.mkd-2}
\end{figure}

\noindent
The jump matrix has the form:
\[
J(k,x,t)=
\begin{cases}
\begin{pmatrix} 
1& \lambda\bar r(k)\ee^{-2\ii\theta}\\
-r(k)\ee^{2\ii\theta}& 1-\lambda |r(k)|^2 
\end{pmatrix}                 &\arg k=0, \ \pi;\\[5mm]
\begin{pmatrix}
1&0\\
c(k)\ee^{2\ii\theta}&1
\end{pmatrix}                 &\arg k=\dfrac{\pi}{3},\  \dfrac{2\pi}{3};\\[5mm]
\begin{pmatrix}
1&\lambda\overline c(\bar k)\ee^{-2\ii\theta}\\
0&1
\end{pmatrix}                 &\arg k=\dfrac{4\pi}{3}, \ \dfrac{5\pi}{3}.
\end{cases}
\]
The proof of equations (\ref{RH}) is a simple algebraic verification
of relations (\ref{IST1})-(\ref{IST2}).

\begin{rem*}
The above Riemann-Hilbert problem is written for the case
when the set of the eigenvalues is empty.
More details about Riemann-Hilbert problem for the mKdV equation and
complete consideration of the initial-boundary value problem
can be found in \cite{BFS}.
\end{rem*}

\begin{rem*}
We consider the problem for the form (\ref{mkd}) of the mKdV
equation in the first quarter $(x\geq 0, t\geq 0)$ of the $xt$-plane.
This problem is different from that studied in \cite{BFS},
which is the initial boundary value problem
for the same equation, but in the second quarter  $(x\leq 0,\ t\geq 0)$
of the $xt$-plane. Scattering (spectral) data for that problem have different
analytic properties. It is well-known
that for the KdV and mKdV equations there are
differences between the initial boundary value problems for $x>0$ and
for $x<0$.
\end{rem*}

Indeed, the kernel of the Marchenko integral equations, which have to
considered now in the domain $-\infty<y<x\le0$, takes the form~:
\begin{align*}	                                         
H(x,t)
&=\frac{1-\lambda}{2}\Bigl(\sum_{k_j\in\Omega_5}
\hat m_j^5 \ee^{\ii k_jx+8\ii k_j^3t}
+\sum_{z_j\in\Omega_4}
m_j^4\ee^{\ii z_jx+8\ii z_j^3t}
+\sum_{z_j\in\Omega_6}
m_j^6\ee^{\ii z_jx+8\ii z_j^3t}\Bigr)
\notag\\
&\quad
+\frac{1}{2\pi}\int_{\partial\Omega_5}
c(k)\ee^{\ii kx+8\ii k^3t}\dd k
+\frac{1}{2\pi}\int_{-\infty}^{\infty}(r(k)+c(k))\ee^{\ii kx+8\ii k^3t}\dd k,
\end{align*}
where $c(k)$ is now meromorphic (analytic for $\lambda=1$) in $\Omega_4\cup\Omega_6$.
For $t=0$ it is easy to prove that
\begin{align*}
&\frac{1}{2\pi}\int_{\partial\Omega_5}
c(k)\ee^{\ii kx}\dd k
=
-\frac{1}{2\pi}\int_{-\infty}^{\infty}c(k)\ee^{\ii kx}\dd k\\
&\qquad\qquad
+\frac{1-\lambda}{2}
\Bigl(-\sum_{z_j\in\Omega_4}
m_j^4\ee^{\ii z_jx}
-\sum_{z_j\in\Omega_6}
m_j^6\ee^{\ii z_jx}
+\sum_{k_j\in\Omega_4}
\hat m_j^4\ee^{\ii k_jx}
+ \sum_{k_j\in\Omega_6}
\hat m_j^6\ee^{\ii k_jx}\Bigr)
\end{align*}
by using analytical properties of function $c(k)$ in the domain
$\Omega_4\cup\Omega_6$ and corresponding residues. Finally we obtain
\[
H(x,0)=\frac{1-\lambda}{2}\left(\sum_{\substack{k_j\in\D{C}_-}}
\hat m_j \ee^{\ii k_jx}\right)
+\frac{1}{2\pi}\int_{-\infty}^{\infty}r(k)\ee^{\ii kx}\dd k.
\]
Marchenko equations with this kernel correspond precisely to the
initial function $u(x)$.

To prove that $q(x,t)$ satisfies the boundary condition we
have also to formulate the Riemann-Hilbert problem for the $t$-equation.
To do so let us define the sectionally meromorphic (analytic for $\lambda=1$)
matrix $N(k,t)$:
\[
N(k,t)=
\begin{cases}
\begin{pmatrix}
\dfrac{\hat\Psi_1^+(t,k)\ee^{-4\ii k^3t}}{p_2^+(k)}&
\hat\varphi_1^-(t,k)\ee^{4\ii k^3t}\\
\dfrac{\hat\Psi_2^+(t,k)\ee^{-4\ii k^3t}}{p_2^+(k)}&
\hat\varphi_2^-(t,k)\ee^{4\ii k^3t}
\end{pmatrix}           &k\in\Omega_1\cup\Omega_3\cup\Omega_5\\[1cm]
\begin{pmatrix}
\hat\varphi_1^+(t,k)\ee^{-4\ii k^3t}&
\dfrac{\hat\Psi_1^-(t,k)\ee^{4\ii k^3t}}{p_1^-(k)}\\
\hat\varphi_2^+(t,k)\ee^{-4\ii k^3t}&
\dfrac{\hat\Psi_2^-(t,k)\ee^{4\ii k^3t}}{p_1^-(k)}
\end{pmatrix}           &k\in\Omega_2\cup\Omega_4\cup\Omega_6,
\end{cases}
\]
where 
$\hat\varphi(t,k)=\begin{pmatrix}\hat\varphi^-(t,k)&\hat\varphi^+(t,k)\end{pmatrix}$ 
and
$\hat\Psi(t,k)=\begin{pmatrix}\hat\Psi^-(t,k)&\hat\Psi^+(t,k)\end{pmatrix}$
are matrix solutions (\ref{hphi}) and (\ref{hatpsi}) of the $t$-equation.
Then it is easy to verify that $N(k,t)$ is a solution of the following
Riemann-Hilbert problem:
\begin{equation} \label{tRH}
N_-(k,t)=N_+(k,t)J^t(k,t), \quad k\in\Sigma
\end{equation}
on the contour $\Sigma$ (Figure 2) oriented as above. The jump matrix
$J^t(k,t)$ has the form:
\[
J^t(k,t)=
\begin{pmatrix}
1&p_-(k)\ee^{4\ii k^3t}\\
-p_+(k)\ee^{-4\ii k^3t}&1-p_-(k)p_+(k)
\end{pmatrix},
\]
where $p_-(k)=p_2^-(k)/p_1^-(k)$, $p_+(k)=p_1^+(k)/p_2^+(k)$, and
$p^\pm_j(k)$ are the entries of the ``scattering'' matrix $P(k)$.
                                                 
We have already proved that
$q(x,0)=u(x)$. The proof that $q(x,t)$ satisfies the boundary values
$q(0,t)=v(t)$, $q_x'(0,t)=v_1(t)$
and $q_{xx}''(0,t)=v_2(t)$ is carried out by using
Riemann-Hilbert problems (\ref{RH}) and (\ref{tRH}) in the same manner as
in \cite{FIS01}.
The main tool in the proof is the existence of an analytic map
from the Riemann-Hilbert problem (\ref{RH}) attached to $M(k,0,t)$ into the
Riemann-Hilbert problem (\ref{tRH}) attached to $N(k,t)$.
Such a  proof for the mKdV equation is precisely given in \cite{BFS}.

\section{Inverse scattering Problem}
\setcounter{equation}{0}

Let $\C{R}$ be scattering data (\ref{calR})
satisfying conditions A\textendash B\textendash C. Then:

\begin{stats*}
\textup{1.} 
The $xt$-integral equation
\begin{align}
& K(x,y,t)+\C{H}(x+y,t)+\int_{x}^{\infty}K(x,z,t)\C{H}(z+y,t)\dd z=0,  \label{xteq} \\
& 0\leq x<y<\infty, \quad 0\leq t<\infty  \notag
\end{align}
with the $2\times 2$ matrix kernel
\[
\C{H}=
\begin{pmatrix}
0 & H(x,t) \\
\lambda {H}(x,t) & 0
\end{pmatrix},
\]
where real scalar function $H(x,t)$ given by \textup{(\ref{H})}, is uniquely
solvable in $L^1(x,\infty)$ for any $x\geq 0$ and 
$t\geq 0$.

\textup{2.} 
The solution $K(x,y,t)$ belongs to 
$C^{\infty}(\D{R}_+\times\D{R}_+\times\D{R}_+)$, it and all its derivatives
decrease faster than any negative power of $x+y$, for $x+y\to\infty
$, and $t$ fixed.

\textup{3.} 
The matrix
\[
\Psi (x,t,k)=
\left[\ee^{-\ii kx\sigma_3}+\int_{x}^{\infty}K(x,y,t)
\ee^{-\ii ky\sigma_3}\dd y\right]\ee^{-4\ii k^{3}t\sigma_3}
\]
satisfies the symmetry conditions
\begin{alignat*}{2}
&\Psi (x,t,k)=\Lambda\bar{\Psi}(x,t,k)\Lambda^{-1}&\quad&\text{for }k\in\D{R}\\
&\Psi^\pm (x,t,k)=\overline{\Psi^\pm(x,t,-\bar k)}&&\text{for }k\in\D{C}_\pm
\end{alignat*}
and is a solution of the $x$-equation $(\ref{xeq})$ with $Q(x,t)$ given by
\begin{equation}
Q(x,t)=\sigma_3K(x,x,t)\sigma_3-K(x,x,t).  \label{Q}
\end{equation}

\textup{4.} 
$\Psi (x,t,k)$ is a  solution of the $x$-and $t$-equations constructed from the matrix
$Q(x,t)$ and its derivative
$Q_x'(x,t)$, $Q_{xx}''(x,t)$, using eqs $(\ref{Q})$, 
$(\ref{xeq})$, $(\ref{teq})$ and $(\ref{PSI})$.

\textup{5.}
The scattering data $\C{R}$ of these compatible
differential equations coincide with the chosen function
$r(k)$, the function $c(k)$ and the numbers
$k_1, k_2, \dots, k_n\in\D{C}_+$, \ $z_1, z_2,\dots,z_m\in\Omega_2$.
\end{stats*}

Statement 1 follows from Lemma \ref{lem.7.1} about the
solvability of the $xt$-integral equations:

\begin{lem}                                             \label{lem.7.1}
Let $\C{R}$ be scattering data satisfying conditions 
\emph{A\textendash B\textendash C}.
Then the $xt$-integral equations
\emph{(\ref{BIE1})-(\ref{H})} have a unique solution 
in $L^1(x,\infty)$.
\end{lem}

\begin{proof}
Under conditions A\textendash B\textendash C the integral
operator of the $xt$-integral equation
is compact in $L^1(x,\infty)$. Then, by Fredholm theory the
$xt$-integral equation has a unique solution if the homogeneous equation
has no non-zero solution.  If a non-zero solution does exists in
$L^1(x,\infty)$, in view of the homogeneity of the integral equation,
it is bounded, hence belongs to $L^2(x,\infty)$.  The integral
operator is clearly skew-Hermitian in $L^2(x,\infty)$, so we obtain a
contradiction, because the only solution in this case is zero.
For $\lambda=1$ the proof is more complicated. For
example, it follows from the solvability of the corresponding
Riemann-Hilbert problem (\ref{RH}). In turn, the unique solvability of the
Riemann-Hilbert problem is proved by the same way as in \cite{XZ}.
\end{proof}

Statement 2 follows from Lemma \ref{lem.7.2}:

\begin{lem}                                            \label{lem.7.2}
Let conditions \emph{A\textendash B\textendash C} be
fulfilled and $(K_1(x,y,t),K_2(x,y,t))$ be the solution of the $xt$-integral
equations \textup{(\ref{BIE1})-(\ref{H})}. 

Then $(K_1(x,y,t),  K_2
(x,y,t))\in C^{\infty}(\D{R}_+ \times\D{R}_+\times\D{R}_+)$.
These functions and all their derivatives decrease
faster than any negative power of $x+y$, for $x+y\to\infty$, 
$t$ fixed. Moreover,
\[
q(x,t)=-2\lambda K_2(x,x,t)
\]

\textup{(a)}
is $C^{\infty}$ in $x$ and $t$,

\textup{(b)}
decreases faster than any negative power of $x$
for $x\to\infty$, $t$ fixed,

\textup{(c)}
is a solution of the mKdV equation with initial function
$q(x,0)=u(x)$ and boundary values $q(0,t)=v(t)$, $q_x'(0,t)=v_1(t)$,
$q_{xx}''(0,t)=v_2(t)$.
\end{lem}

\begin{proof}
According to Lemma \ref{lem.7.1} the $xt$-integral equations have a solution
\[
K(x,y,t)=\begin{pmatrix}
K_1(x,y,t)&\lambda K_2(x,y,t)\\
K_2(x,y,t)&K_1(x,y,t)
\end{pmatrix}
\]
which
belongs to $L^1(x,\infty)$. By condition A, the kernel
$H(x+y,t)$ is in $C^{\infty}(\D{R}_+\times\D{R}_+)$
and is fast decreasing as $x+y\to\infty$
(end of Section 4). Therefore $(K_1(x,y,t),K_2(x,y,t))$ is in
$C^{\infty}(\D{R}_+\times\D{R}_+\times\D{R}_+)$ and
vanishes faster than any negative power of $x+y$
as $x+y\to\infty$, $t$ fixed. The same is true for all
their derivatives and for $q(x,t)$.
It is clear that $H(x,t)$ satisfies:
\[
\frac{\partial}{\partial t} H(x,t) + 8\frac{\partial^3}{\partial
x^3}H(x,t)=0.
\]
Then it is well-known that $q(x,t)$ solves the mKdV
equation, cf.\  \cite{AS,ZSh}.
The fact that $q(x,t)$
satisfies the boundary conditions was discussed at the end of Section 5.
\end{proof}

\begin{proof}[Proof of Statements $3$ and $4$]
The proof of Statement 3 is well-known  \cite{AS,FT}.
In particular, formula (\ref{Q}) follows from equation
(\ref{sigma}). Statement 4
is also true. Indeed, due to Lemma \ref{lem.7.2},
the function $q(x,t)$ is a solution of the mKdV equation (\ref{mkd}). Hence the
constructed $x$- and
$t$-equations are compatible. 
Therefore, according to Lemma \ref{lem.2.1}, the
matrix-valued function $\Psi(x,t,k)$ solves the $t$-equation.
\end{proof}

\begin{proof}[Proof of Statement $5$]

Let us consider compatible solutions  $\Psi(x,t,k)$ (\ref{PSI}),
$\Phi(x,t,k)$ (\ref{PHI}) and $Y(x,t,k)$ (\ref{Y})
of the constructed $x$- and $t$-equations.
Let $\tilde{\C{R}}$
be the corresponding scattering data.
We have to show that $\tilde{\C{R}}\equiv\C{R}$.

First of all one can find that the matrix
$K(x,y,0)$ solves the $xt$-integral equations ($t=0$) with
kernel $\C{H}(x,0)$ generated by
\[
\tilde H(x)=\sum_{\substack{\tilde k_j\in{\D{C}_+}\\
                  \tilde s_2^+(\tilde k_j)=0}}
   \tilde m_j\ee^{\ii\tilde k_j x}
   +\frac{1}{2\pi}\int_{-\infty}^{\infty}\tilde r(k)\ee^{\ii kx}\dd k.
\]
On the other hand the matrix $K(x,y,0)$ solves the same
integral equation with kernel generated by
\[
H(x)=\sum_{\substack{k_j\in{\D{C}_+}\\
                     s_2^+(k_j)=0}}
          m_j\ee^{\ii k_j x}
+\frac{1}{2\pi}\int_{-\infty}^{\infty}r(k)\ee^{\ii kx}\dd k,
\]
since
\[
\frac{1}{2\pi}\int_{\partial\Omega_2}
c(k)\ee^{\ii kx}\dd k=
    \sum_{\substack{k_j\in\Omega_2\\
                    s_2^+(k_j)=0}}m_j^2\ee^{\ii k_jx}
    -\sum_{\substack{z_j\in\Omega_2\\
                     r_1^-(z_j)=0}}m_j^2\ee^{\ii z_jx}.
\]
Hence $F(x)=\tilde H(x)-H(x)$ is a solution
of the homogeneous Volterra integral equation
\[
F(2x)+\int_{2x}^{\infty}K_1(x,y-x,0)F(y)\dd y=0,
\]
which yields the identity $F(x)\equiv 0$. Therefore
$\tilde k_j=k_j$,  $\tilde m_j=m_j$ (hence $\tilde m^1_j=m^1_j$  
$\tilde m^3_j=m^3_j$)  and
$\tilde r(k)\equiv r(k)$  for $k\in\D{R}$.
For $t>0$ by using (\ref{H}) we shall obtain now the relation
$ F(x,t)=\tilde H(x,t)-H(x,t)\equiv0$, that yields
$\tilde c(k)\equiv c(k)$ and $\tilde z_j =z_j$,  $\tilde m^2_j=m^2_j$.
Hence $\tilde{\C{R}}\equiv\C{R}$.
\end{proof}

All statements on the
inverse scattering problem are proved.
The main theorem is proved.


\end{document}